\documentclass[11pt, letterpaper]{article}

\usepackage{adjustbox}
\usepackage{algorithm}
\usepackage{algpseudocode}
\usepackage{amsmath}
\usepackage{amssymb}
\usepackage{amsthm}
\usepackage{authblk}
\usepackage{booktabs}
\usepackage{cancel}
\usepackage{caption}
\usepackage{comment}
\usepackage{diagbox}
\usepackage{dsfont}
\usepackage{eurosym}
\usepackage[OT1]{fontenc}
\usepackage{float}

\usepackage[top=1.75cm, bottom=3.75cm, left=2cm, right=2cm]{geometry}
\usepackage{natbib}
\usepackage[colorlinks=true, allcolors=blue]{hyperref}

\usepackage{mathrsfs}
\usepackage{mathtools}
\usepackage{multirow}
\usepackage{nicematrix}
\usepackage{siunitx}
\usepackage{subcaption}
\usepackage{upquote}
\usepackage{xcolor}

\numberwithin{equation}{section}

\newtheorem{lemma}{Lemma}[section]

\newtheorem{remark}{Remark}[section]
\newtheorem{theorem}{Theorem}[section]
\newtheorem{corollary}{Corollary}[section]
\theoremstyle{remark}

\DeclareMathOperator{\N}{\mathbb{N}}
\DeclareMathOperator{\R}{\mathbb{R}}
\DeclareMathOperator{\Z}{\mathbb{Z}}

\DeclareMathOperator{\Corr}{\mathrm{Corr}}
\DeclareMathOperator{\Erw}{\mathds E}
\DeclareMathOperator{\Prob}{\mathds P}

\DeclareMathOperator{\tr}{\mathrm{tr}}

\begin{document}

\title{Pivotal inference for linear predictions in stationary processes}

\author{
Holger Dette\thanks{Faculty of Mathematics, Ruhr-Universit\"at Bochum, \href{mailto:holger.dette@rub.de}{holger.dette@rub.de}} 
\qquad
Sebastian Kühnert\thanks{Faculty of Mathematics, Ruhr-Universit\"at Bochum, \href{mailto:sebastian.kuehnert@rub.de}{sebastian.kuehnert@rub.de}}
}

\date{November 13, 2026}

\maketitle

\begin{abstract} 
In this paper we develop pivotal inference for the final (FPE) and relative final prediction error (RFPE) of linear forecasts in stationary processes. Our approach is based on a self-normalizing technique and avoids the estimation of the asymptotic variances of the empirical autocovariances. We provide pivotal confidence intervals for the (R)FPE, develop estimates for the minimal order of a linear prediction that is required to obtain a prespecified forecasting accuracy and also propose (pivotal) statistical tests for the hypotheses that the (R)FPE exceeds a given threshold. Additionally, we provide pivotal uncertainty quantification for the commonly used coefficient of determination $R^2$ obtained from a linear prediction based on the past $p \geq 1$ observations and develop new (pivotal) inference  tools for the partial autocorrelation, which do not require the assumption of an autoregressive process.
\end{abstract} 



\section{Introduction}

Linear forecasting is an important and central  technique of time series analysis due to its simplicity, ease of interpretation, and well-established theoretical properties. In the simplest case, it aims  for  predicting future values, such as $X_n$ of a stationary temporal process $(X_k)_{k \in \mathbb{Z}} $ using a linear combination, say $\widehat X_{n,p},$ of its past observations $X_{n-1}, \ldots , X_{n-p}$. An important problem is the estimation of the prediction error $\Erw\! |X_n- \widehat X_{n,p}|^2$ as it enables the construction of prediction intervals. Prediction error estimates are also often used to obtain reasonable models for fitting the data, which corresponds to the  choice of the appropriate order $p$ for linear predictions. Numerous other criteria for model selection have been proposed in the literature. Prominent examples include the Akaike Information Criterion (AIC) \citep{Akaike1974}, which aims to minimize the mean squared final prediction error (FPE) \citep{Akaike1969}, the Bayesian Information Criterion (BIC) \citep{schwarz}, and the criterion of \cite{HannanQuinn1979}. Subsequent contributions have further advanced this area; see, e.g., the comprehensive treatment in \cite{claeskensHjort08}. Most of these works aim to establish consistency in selecting the order of potentially infinite-order autoregressive models.

In this paper, we adopt a different perspective on assessing the quality of linear prediction in a centered, stationary linear process $(X_k)_{k \in \mathbb{Z}}$, focusing on the {\it minimum distance}
\begin{align}\label{Mp}
    M_p \coloneqq \min_{\xi_1, \dots, \xi_p \in \mathbb{R}}\, \Erw\!\Big(X_n - \sum_{i=1}^p \xi_i X_{n-i}\Big)^2 = \min_{\xi_1, \dots, \xi_p \in \mathbb{R}}\, \Erw\!\Big(X_p - \sum_{i=1}^p \xi_i X_{p-i}\Big)^2,
\end{align}
with $M_0 \coloneqq \Erw(X_0^2),$ measuring the mean squared deviation between $X_n$ and its optimal linear predictor based on the past $p \geq 1$ observations. Note that $M_p$ is the population counterpart of the mean squared final prediction error considered in \cite{Akaike1969}. We will develop pivotal statistical inference for $M_p$ and use it to quantify uncertainty in estimates of related measures commonly used to assess the quality of linear predictions. Specifically, we consider the measures
\begin{align}\label{det91}
    S_p  = {M_p \over M_0} \qquad\quad\text{and}\qquad\quad
    Q_p  = {M_p \over M_{p-1}}, 
\end{align}
which are normalized versions of \eqref{Mp}. Both measures have been extensively studied and applied by various authors; see, for instance, \cite{BrockwellDavis1991, Hannan1970, Ramsay1974}. The measure $S_p$ considers the mean squared error of the best linear prediction relative to the magnitude of $\Erw(X_n^2)$ and is the population version of the {\it relative final prediction error} (RFPE) \citep[see][]{Akaike1969}. The measure $S_p$ is directly related to the \emph{coefficient of determination} $R^2_p$ 
through $S_p=1-R^2_p,$ and $Q_p$ fulfills $Q_p = 1 - \kappa_p^2$, where $\kappa_p$ denotes the $p$th partial autocorrelation. We will develop several pivotal inference tools for these quantities, which are briefly described for $S_p$ here. First we are interested in (pivotal) confidence intervals for $S_p$ and $R_p^2$. For example, we will provide  pivotal uncertainty quantification for the commonly used coefficient of determination $R^2_p$ obtained from a linear prediction based on the past $p$ observations. Second, we construct a test for the hypothesis that $S_p$ is sufficiently small, that is 
\begin{align} \label{revdet1}
    H_0 : S_p \geq \Delta \quad \text{ vs. } \quad H_1 : S_p < \Delta,
\end{align}
where $\Delta\in(0,1)$ is a given threshold. Rejecting the null hypothesis in favor of $S_p<\Delta$ implies a coefficient of determination exceeding $1-\Delta$, with the type I error controlled by a pivotal distribution. Third, we aim to estimate, with a statistical guarantee, the minimal lag at which the coefficient of determination reaches a prespecified value $\nu$, that is
\begin{equation}\label{2.18}
    p^* = \min \big\{ p \in \mathbb{N} \;\big|\; R_p^2 \geq \nu  \big\} = \min \big\{ p \in \mathbb{N} \;\big|\; S_p \leq 1 - \nu \big\}. 
\end{equation}
Finally, we assess the adequacy of a linear predictor of order $p_0$ by formally testing the hypotheses 
\begin{align}\label{eq_Hypos_Intro}
    H_0 : p_0 < p^*  \quad\text{ vs. }\quad H_1 : p_0 \geq p^*.    
\end{align}
These objectives can be addressed by expressing $S_p$ as a function of the autocovariances $\gamma_0,\ldots,\gamma_p$, say $S_p=f(\gamma_0,\ldots,\gamma_p)$, and estimating it by $\widehat S_p=f(\widehat\gamma_0,\ldots,\widehat\gamma_p)$ using sample autocovariances. The joint asymptotic normality of $\sqrt{N}(\widehat\gamma_0-\gamma_0,\ldots,\widehat\gamma_p-\gamma_p)^\top$, where $N$ is the sample size, then yields the asymptotic normality of $\widehat S_p$ via the (functional) delta method \citep[see][]{vanderVaart1998}. For general linear processes, however, the asymptotic covariance matrix depends on the entire autocovariance sequence $(\gamma_k)_{k\in\mathbb{N}_0}$ \citep[see, e.g.,][Chapter 7]{BrockwellDavis1991}. Consequently, the asymptotic variance of $\sqrt{N}(\widehat S_p-S_p)$ is difficult to estimate and typically requires regularization.

To circumvent these problems, we consider a different approach based on self-normalization, allowing for pivotal inference regarding the measures \eqref{Mp}--\eqref{det91}. While self-normalization has received considerable attention for testing whether a parameter vanishes \citep[see, e.g.,][]{lobato,shaozhang2010,shao2015}, it has been much less explored for composite hypotheses of the form \eqref{revdet1} and estimation and testing problems of the types \eqref{2.18} and \eqref{eq_Hypos_Intro}, respectively. The paper is organized as follows. Section~\ref{sec2} introduces the general model and the foundations of our method. Section~\ref{sec3} develops statistical inference for $S_p$ and $R_p^2$, while Section~\ref{sec4} presents corresponding results for $Q_p$ and extends existing inference methods for partial autocorrelations $\kappa_h$. Section~\ref{sec5} outlines extensions to the multivariate setting. All proofs are deferred to the appendix.

\section{Sequential estimation of the  final prediction error} \label{sec2}

This section introduces a sequential estimator of the measure $M_p$ in \eqref{Mp}. Its properties are crucial for developing inference tools for the associated normalized measures $S_p$ and $Q_p$ in \eqref{det91}, considered in Sections~\ref{sec3}--\ref{sec4}. Although $M_p$ is not the main object of interest in this paper, we illustrate the self-normalizing approach for $M_p$, as it is technically more transparent than for $S_p$ and $Q_p.$ 

Let $(X_k)_{k \in \mathbb{Z}} $  denote a centered  stationary process with finite variance. Motivated by Wold's decomposition theorem \citep[see][Section 5.7]{BrockwellDavis1991}, we assume that
\begin{align}\label{eq:(X_k) lin process}
    X_k = \sum^\infty_{j=0}\,\theta_{\!j}\varepsilon_{k-j},\quad k\in\Z,  
\end{align}
where the coefficients $\theta_{\!j}$ satisfy  $0 < \sum^\infty_{j=1} j| \theta_j| < \infty$ and $(\varepsilon_k)_{k\in\Z}$ is an i.i.d.\ process of innovations satisfying $\Erw(\varepsilon^4_0)<\infty$, $\Erw(\varepsilon_0)=0$ and $\Erw(\varepsilon_0^2)=1$. We denote by $\gamma_h = \mathbb{E}(X_0 X_h)$  
the {\it autocovariance} at lag $h \in \mathbb{N}_0$ of such a process, and by $\kappa_h$ the {\it partial autocorrelation} at lag $h \in \mathbb{N},$ which is defined by $\kappa_1 \coloneqq \Corr(X_0,X_1),$ and $\kappa_h \coloneqq \Corr(X_0 - \widehat{X}_0, X_h - \widehat{X}_h)$ for $h>1,$ where $\widehat{X}_0$ and $\widehat{X}_h$ denote the linear projections of $X_0$ and $X_h$ onto $\{X_1, X_2, \dots, X_{h-1}\}$, respectively. Further, we interpret $\langle X_0, X_h \rangle = \Erw (X_0X_h )$ as the common inner product on the Hilbert space ${\cal H} $ of all centered square integrable random variables. Classical results from approximation theory 
\citep[see][p.\ 16]{Achieser1956} provide a simple representation for the final population prediction error as a ratio of two \emph{Gram determinants}, that is
\begin{align}\label{Def M_p}
    M_p = \underset{\xi_1, \dots, \xi_p\,\in\R}{\min}\,\Erw\!\Big(X_p - \sum^p_{i=1}\,\xi_iX_{p-i}\Big)^2 ~=~
    \frac{\det(G_p)}{\det(G_{p-1})}\,,
\end{align}
where $G_p \coloneqq (\gamma_{j-i})_{i,j=0}^p$ is the Toeplitz matrix of the autocovariances $\gamma_0, \dots, \gamma_p$, and $\det(G_{-1})=1$ by convention. Throughout this paper we assume that all matrices $G_p$ are non-singular, which, e.g., is satisfied if $\gamma_h \to 0$ as $h \to \infty $ \citep[][Proposition~5.1.1]{BrockwellDavis1991}, thus $M_p > 0.$ Next, let
\begin{align}\label{Def sample autocov function}
    \widehat{\gamma}_h \coloneqq \frac{1}{N}\sum_{i=1}^{N-|h|} X_i X_{i+|h|}, \qquad |h| < N,
\end{align}
denote the usual estimator of the autocovariance $\gamma_h$ from a sample $X_1, \ldots ,X_N.$ A  canonical estimator for $M_p$ is then given by 
\begin{align}\label{Estimate for M_p}
    \widehat{M}_p \coloneqq \frac{\det(\widehat{G}_p)}{\det(\widehat{G}_{p-1})}\,,
\end{align}
with Toeplitz matrices $\widehat{G}_p\coloneqq (\widehat{\gamma}_{j-i})^p_{i,j=0}$ which are non-negative definite, and positive definite (denoted $\succ 0$) whenever $\widehat{\gamma}_0 > 0$ \citep[cf.][Section~7.2]{BrockwellDavis1991}. The asymptotic distribution of $\widehat M_p$ can then be derived from that of the vector $\widehat{\boldsymbol{\gamma}}_{\!p} \coloneqq (\widehat{\gamma}_0, \dots, \widehat{\gamma}_p)^{\!\top}\!,$ since, with $\boldsymbol{\gamma}_p \coloneqq (\gamma_0, \dots, \gamma_p)^{\!\top}\!$, we have
\begin{align}
    \sqrt{N}(\widehat{\boldsymbol{\gamma}}_{\!p} - \boldsymbol{\gamma}_p) &\stackrel{d}{\longrightarrow} \mathcal{N}(0, \Sigma_p),\label{AN for vector of autocov}
\end{align}
where the elements of the covariance matrix $\Sigma_p = (\Sigma_{p,ij})^p_{i,j=0}\in\R^{(p+1)\times (p+1)}$ are given by
\begin{align}   
    \Sigma_{p,ij} &\coloneqq  \gamma_i \gamma_j \Erw(\varepsilon^4_0 - 3) + \sum_{k=-\infty}^\infty \gamma_k \gamma_{k-i+j} + \gamma_{k+j} \gamma_{k-i}\,\label{Def of Sigma_ij}
\end{align}
\citep[see][Chapter~7]{BrockwellDavis1991}. A direct application of the delta method yields
\begin{align}\label{det95}
    \sqrt{N}\big(\widehat{M}_p-M_p\big)\, \stackrel{d}{\longrightarrow} \,\mathcal{N}(0,\tau^2_p),
\end{align}
where $\tau^2_p \coloneqq (\nabla M_{p,\boldsymbol{\gamma}_p})^\top \Sigma_p\,\nabla M_{p,\boldsymbol{\gamma}_p}$, and $\nabla M_{p,\boldsymbol{\gamma}_p}$ is the gradient of $M_p$ evaluated at $\boldsymbol{\gamma}_p$ (gradients are throughout taken as column vectors). Valid inference requires a consistent estimate of the asymptotic variance $\tau^2_p$, a challenging task since it depends on the full dynamics of $(X_k)$ through the autocovariances $(\gamma_h)_{h\in \mathbb{N}_0}$. While estimation of $\nabla M_{p,\boldsymbol{\gamma}_p}$ is straightforward, estimating $\Sigma_p$ in \eqref{AN for vector of autocov} is more delicate. Standard approaches truncate the series in \eqref{Def of Sigma_ij} at some $k_n \in \mathbb{N}$ \citep[see, e.g.,][]{lee2003}, replacing the unknown autocovariances and fourth moment of the innovations with corresponding estimates. However, choosing $k_n$ and reliably estimating $\Erw(\varepsilon_0^4)$ remain non-trivial.

To circumvent these difficulties, we pursue an alternative route and derive a pivotal limiting distribution for $\widehat M_p - M_p$ after normalization by a suitable factor. For this purpose, define
\begin{align}\label{Estimate gamma_h uni-variate}
    \widehat{\gamma}_h(\lambda)\coloneqq \frac{1}{N}\!\sum^{\lfloor \lambda(N-|h|)\rfloor}_{i=1}\!X_iX_{i+|h|}, \quad \lambda \in [0,1],
\end{align}
as a sequential estimator of $\lambda\gamma_h$ ($|h| < N$), where $\widehat{\gamma}_h(1)$ equals $\widehat \gamma_h$ in \eqref{Def sample autocov function}. Moreover, set
\begin{align} \label{hdrev1}
    \widehat{V}_{M_p} \coloneqq \int_0^1 \big|\widehat{M}_p(\lambda) - \lambda \widehat{M}_p\big|\,\mathrm{d}\lambda,
\end{align}
with
\begin{align}\label{Def Mp gamma}
    \widehat M_p(\lambda)
    \coloneqq
    \begin{cases}
        \,\dfrac{\det(\widehat G_p(\lambda))}
        {\det(\widehat G_{p-1}(\lambda))},
        & \widehat G_p(\lambda)\succ0,\\[2mm]
        \,0, & \text{otherwise},
    \end{cases}
    \qquad \lambda\in[0,1],
\end{align}
where
$\widehat G_p(\lambda)
\coloneqq(\widehat\gamma_{j-i}(\lambda))_{i,j=0}^p.$
In particular, $\widehat M_p(1)$ coincides with
$\widehat M_p$ in \eqref{Estimate for M_p}
whenever $\widehat G_p(1)\succ0$.

Our first main result establishes the weak convergence of the process 
\begin{align}\label{det94}
    \big\{\pmb{\mathcal{I}}_N(\lambda) \big\}_{\!\lambda\in[0,1]} \coloneqq 
    \sqrt{N}\big\{\big(\widehat{M}_0(\lambda),  \dots, \widehat{M}_p(\lambda)\big)^{\!\top} - \lambda \big(M_0, \dots, M_p\big)^{\!\top}\big\}_{\lambda\in[0,1]}\,,
\end{align}
which is fundamental for developing pivotal inference for $S_p$ and $Q_p$ in \eqref{det91}. For a precise statement, let $\ell^\infty([0,1])$ be the space of bounded functions $f\colon[0,1]\to\R$ with norm $\|f\|_\infty=\sup_{\lambda\in[0,1]}|f(\lambda)|$, and let
$$
    \ell^{\infty, p+1}([0,1]) \coloneqq \left\{\boldsymbol{f}_{\!p} = (f_0, \dots, f_p)^{\!\top} : [0,1] \to \R^{p+1} ~\big|~ f_j \in \ell^{\infty}([0,1]),~ j=0,\dots,p \right\}
$$
denote the space of bounded functions $\boldsymbol{f}_{\!p}\colon [0,1]\to \R^{p+1}$ with norm $\|\boldsymbol{f}_{\!p}\|_\infty=\sup_{\lambda\in[0,1]}\sup_{0\leq i\leq p}|f_i(\lambda)|$. Throughout, $\rightsquigarrow$ denotes weak convergence in $\ell^{\infty}([0,1])$ or $\ell^{\infty, p+1}([0,1])$, with the underlying space being clear from the context \citep[see][for convergence of processes]{VaartWellner2023}.

\begin{theorem}\label{M0, ..., Mp cvgc} The process in  \eqref{det94} satisfies
\begin{align*}
    \big\{\pmb{\mathcal{I}}_N(\lambda) \big\}_{\!\lambda\in[0,1]}
    \rightsquigarrow \pmb{\mathcal{I}}\coloneqq \big\{\mathcal{M}_{p,\boldsymbol{\gamma}_p}\,\Sigma_p^{1/2}\,\pmb{\mathbb{B}}(\lambda)\big\}_{\!\lambda\in[0,1]}\,,  
\end{align*}
where $\Sigma_p\in\R^{(p+1)\times(p+1)}$ and $\mathcal{M}_{p,\boldsymbol{\gamma}_p}\in\R^{(p+1)\times(p+1)}$ are the matrices in \eqref{AN for vector of autocov} and \eqref{kue7}, respectively, and $\pmb{\mathbb{B}}(\lambda)\coloneqq(\mathbb{B}_0(\lambda),\dots,\mathbb{B}_p(\lambda))^{\!\top}$ is a vector of independent standard Brownian motions on $[0,1]$.
\end{theorem}

The usefulness of Theorem~\ref{M0, ..., Mp cvgc} is illustrated by the following corollary, which directly follows from the continuous mapping theorem.
\begin{corollary}\label{piv_M_cor}  
Let $\Sigma_p$ in  \eqref{AN for vector of autocov} be non-singular. Then, 
\begin{align}\label{def W}
    \frac{\widehat{M}_p-M_p}{\widehat{V}_{M_p}}\;\stackrel{d}{\longrightarrow} \; W \coloneqq  
    \frac{\mathbb{B}(1)}{\int^1_0 |\,\mathbb{B}(\lambda) - \lambda\mathbb{B}(1)|\,\mathrm{d}\lambda}\,,     
\end{align}
where $\mathbb{B}=\{\mathbb{B}(\lambda ) \}_{\lambda \in [0,1]} $
is a standard Brownian motion on the interval $[0,1].$
\end{corollary}

\noindent The denominator in the limiting distribution is almost surely positive. Using the Karhunen--Loève expansion for the Brownian motion it can be shown  that the distribution of $W$ is symmetric. Thus,
$$
\left[ 
    \widehat{M_p} -  q_{1-\alpha/2}(W) {\widehat{V}_{M_p}}~,~ \widehat{M_p} +  q_{1-\alpha/2}(W) {\widehat{V}_{M_p}}
\right]
$$
defines a pivotal asymptotic confidence interval for the measure $M_p,$ where $q_{1-\alpha/2}(W)$ refers to the $(1-\alpha/2)$-quantile of the distribution of $W.$ 

\begin{remark}\label{rem:L^p-m}
~~
    {\rm 
    \begin{itemize}
        \item[(a)]The linear representation \eqref{eq:(X_k) lin process} assumes i.i.d.\ innovations, a standard condition in sequential estimation of autocovariances \citep[see, e.g.,][]{lee2003,BERKES20092044}. The results, however, remain valid under weaker assumptions. In particular, they hold for processes with white noise innovations provided $(X_k)$ is $L^4$-$m$-approximable \citep{HoermannKokoszka2010}, a condition satisfied for independent innovations with finite fourth moments if $\sum_{j=1}^\infty j|\theta_j| < \infty$.
        
        \item[(b)]Self-normalization is a common tool for avoiding the estimation of nuisance parameters in statistical inference \citep[see][for early references]{lobato,shaozhang2010,shao2015}. The self-normalizing statistic in \eqref{hdrev1} differs from $\widetilde V_{M_p} = ( \int_0^1 |\widehat{M}_p(\lambda) - \lambda \widehat{M}_p|^2\,\mathrm{d}\lambda)^{1/2}$, the analog of the statistics used in these references. The proofs show that similar results can be obtained if $\widehat V_{M_p}$ in \eqref{def W} is replaced by $\widetilde V_{M_p}$. Moreover, time-symmetric self-normalization methods as proposed by \cite{LavitasZhang2018} could also be used. Most existing work develops self-normalization techniques to obtain an asymptotically pivotal distribution of $\widehat M_p$ when $M_p=0$. In contrast, the self-normalization in \eqref{def W} addresses the case $M_p>0$. This requires a different asymptotic analysis of the statistic in \eqref{def W}, since one cannot work under the null hypothesis $M_p=0$; see the appendix for details. Related methods were recently applied by \cite{Detteetal2020} to functional data and by \cite{DelftDette2024} to spectral analysis of non-stationary data.
        
        \item[(c)]For technical convenience, the autocovariance estimators in \eqref{Def sample autocov function} and their sequential counterparts $\widehat{\gamma}_h(\lambda)$ in \eqref{Estimate gamma_h uni-variate} assume a centered process $(X_k)_{k\in\mathbb{Z}}$. However, all results remain valid for the estimators
        \begin{align*}
            \widehat{\gamma}_h \coloneqq
            \frac{1}{N}\sum_{i=1}^{N-|h|}
            (X_i-\bar{X}_N)(X_{i+|h|}-\bar{X}_N),
            \qquad |h|<N,
        \end{align*}
        and their sequential counterparts, where $\bar{X}_N = \frac 1N\sum_{i=1}^N X_i$ denotes the sample mean.
    \end{itemize}
    }
\end{remark}

\section{The measure $S_p$ and the coefficient of determination $R^2_p$}\label{sec3}    

This section derives several inference tools for the relative final prediction error in \eqref{det91},
\begin{align} \label{det96}
    S_p \coloneqq \frac{M_p}{M_0} = \frac{M_p}{\gamma_0} = 1 - R_p^2,
\end{align}
where $R_p^2$ is the \emph{coefficient of determination}. We begin with an analogue of Corollary \ref{piv_M_cor}. To this end, we define the sequential estimator
\begin{align}\label{det80}
    \widehat{S}_p(\lambda)\coloneqq
    \begin{cases}
    \,\displaystyle
    \frac{\widehat{M}_p(\lambda)}
    {\widehat{M}_0(\lambda)},
    & \lambda\in(0,1],\;
    \,\widehat{M}_0(\lambda)>0,\\[6pt]
    \,0,
    & \lambda\in(0,1],\;
    \,\widehat{M}_0(\lambda)=0,\\
    \,1,
    & \lambda=0,
    \end{cases}
\end{align}
where $\widehat{M}_p(\lambda)$ is defined
in \eqref{Def Mp gamma}, and we write
$\widehat{S}_p\coloneqq\widehat{S}_p(1)$. We then consider the self-normalizer
\begin{align}
    \widehat{V}_{S_p} \coloneqq \int_0^1 \lambda\big|\widehat{S}_p(\lambda) - \widehat{S}_p\big|\,\mathrm{d}\lambda. \label{Vhat_S}
\end{align}
\begin{theorem}\label{thm: Distr. free limit of M'}
Suppose $\Sigma_p$ in \eqref{AN for vector of autocov} is non-singular and $R_p^2>0$. Then, with $W$ as in \eqref{def W},
\begin{align}\label{kue11a}
    \frac{\widehat{S}_p-S_p}{\widehat{V}_{S_p}}
    \stackrel{d}{\longrightarrow} W.
\end{align}
\end{theorem}

\subsection{Confidence intervals and  testing  relevant  hypotheses} 
\label{sec31}

A pivotal asymptotic confidence interval for the relative final prediction error $S_p>0$ is given by
\begin{align} \label{det97}
\Big [ 
    \widehat S_p  -  q_{1-\alpha/2}(W){\widehat{V}_{S_p}}~,~ \widehat S_p  +  q_{1-\alpha/2}(W){\widehat{V}_{S_p}}
    \Big ],
\end{align}
which directly extends to a pivotal confidence interval for the coefficient of determination $R^2_p$, namely
$$
    \Big[\widehat R_p^2  -  q_{1-\alpha/2}(W) {\widehat{V}_{S_p}}~,~ \widehat R_p^2  +  q_{1-\alpha/2}(W) {\widehat{V}_{S_p}}\Big].
$$
In other words, our approach provides pivotal uncertainty quantification for the commonly used $R^2$ in linear prediction based on the past $p \geq 1$ observations. Next we construct a test for the hypotheses 
\begin{align}
\label{det103}
H_0 : S_p \geq \Delta \quad\text{ vs. }\quad H_1 : S_p < \Delta, 
\end{align}
or equivalently 
\begin{align}
    \label{det103a}
H_0 : R_p^2 \leq 1 - \Delta \quad\text{ vs. }\quad H_1 : R_p^2 >  1- \Delta,
\end{align}
where $\Delta\in(0,1)$ is a prespecified threshold. To test whether $R_p^2$ exceeds $1-\Delta$ while controlling the type~I error, we reject the null hypothesis in \eqref{det103} or \eqref{det103a} whenever
\begin{align} 
\label{det102}
    \widehat{S}_p \leq \Delta +  q_{\alpha}(W)\widehat{V}_{S_p}.
\end{align}
The next result establishes that this procedure yields a consistent asymptotic level $\alpha$-test.

\begin{corollary}\label{Cor:Asymptotic level test - distr. free} Under the assumptions of Theorem \ref{thm: Distr. free limit of M'}, we have 
\begin{align*}
    \lim_{N\rightarrow\infty}
    \Prob\!\big(\widehat{S}_p \leq \Delta + q_{\alpha }(W)\widehat{V}_{S_p} \big) =     
    \begin{cases}
        \;1,      & \text{if } S_p < \Delta,\\ 
        \;\alpha, & \text{if } S_p = \Delta,\\
        \;0,      & \text{if } S_p > \Delta.
    \end{cases}
\end{align*}
\end{corollary}

\begin{remark}\label{remchoceofdelta} 
{\rm Testing hypotheses of the form \eqref{det103} or \eqref{det103a} requires specifying the threshold $\Delta$, which is application-specific and should be carefully justified. For example, to assess whether a linear predictor of order $p$ attains a coefficient of determination exceeding 80\%, a natural choice is $\Delta=0.2$. If the scientific context does not provide a natural choice of the relevance threshold $\Delta$, the nested structure of the hypotheses in \eqref{det103} can be exploited. In particular, rejection of the null hypothesis by the test \eqref{det102} at $\Delta=\Delta_0$ also implies rejection for all $\Delta \geq \Delta_0$. By the sequential rejection principle  we can  simultaneously    test  the hypothesis  \eqref{det103} and determine the minimal $\Delta$ for which the null hypothesis is rejected, that is 
\begin{align*}
    \widehat \Delta_\alpha\coloneqq \!
        \,\inf\!\big\{\Delta \in [0,1]\,\big|\,
        \widehat{S}_p\leq\Delta+q_\alpha(W)\widehat{V}_{S_p}\big\}
\end{align*}
where we use the convention $\inf\emptyset=\infty$.
Since rejection of the null hypothesis occurs
precisely when $\Delta\geq\widehat\Delta_\alpha$
for $\Delta\in[0,1]$, the threshold
$\widehat\Delta_\alpha$ quantifies the evidence
against the null hypothesis, with smaller
values indicating stronger evidence that
the final prediction error is small.
}
\end{remark}

\subsection{Estimating the order for linear predictions}
\label{sec32} 

Recall that $p^*$ in \eqref{2.18} denotes the smallest order in linear prediction for which the coefficient of determination $R_p^2$ reaches a threshold $\nu$ (equivalently, the relative final prediction error $S_p$ is at most $1-\nu$). As $S_p = \prod_{h=1}^p (1-\kappa_h^2)$ \citep[cf.][p.\ 22]{Hannan1970}, where $\kappa_h$ is the $h$th partial autocorrelation, it follows that
$$
    \lim_{p \to \infty} R^2_{p} = 1- \lim_{p \to \infty} S_{p}  = R_\infty^2\coloneqq 
    1 -\prod_{h=1}^\infty  (1- \kappa_h^2),
$$ 
and $R_\infty^2=1$ if and only if $\sum_{h=1}^\infty \kappa^2_h=\infty$. Hence, $p^*$ in \eqref{2.18} is well-defined for all $\nu \in (0,R_\infty^2)$ and throughout this section we only consider this case. We define the estimator
\begin{equation}\label{hx3}
  \widehat{p} = \min \big\{p\in\{1,\dots,N-1\} ~\big|~ \widehat{S}_p < 1-\nu-q_\alpha(W)\widehat V_{S_p} \big\},
\end{equation}
where we use the $\min\emptyset=N$ here, and where $\widehat{S}_p=\widehat{S}_p(1)$ and $\widehat V_{S_p}$ are given in \eqref{det80} and \eqref{Vhat_S}, respectively. Moreover, this estimator satisfies
\begin{theorem}
    \label{thm3} 
Under the assumptions of Theorem \ref{thm: Distr. free limit of M'},
\begin{eqnarray*}
&&  \lim_{N \to \infty}  \mathbb{P} \big (\widehat{p} < p^* \big ) = 0  \quad\text{ and  }\quad \limsup_{N \to \infty}\, \mathbb{P} \big (\widehat{p} > p^* \big ) \leq  \alpha .
\end{eqnarray*}
Moreover, if $\alpha=\alpha_N$ in \eqref{hx3} depends on the sample size $N$ with $\alpha_N \to 0$ and $|q_{\alpha_N}(W)| = o_{\Prob}(\sqrt{N})$, then
$$
\lim_{N \to \infty}  \mathbb{P} \big (\widehat{p} \not =  p^* \big )  = 0\,.
$$
\end{theorem} 

\begin{remark}\label{precision}
    {\rm In finite samples, identifying $p^*$ becomes more difficult as $d^*\coloneqq\max\{S_{p^*-1}-S_{p^*}, S_{p^*}-S_{p^*+1}\}$ decreases. As the proof of Theorem \ref{thm3} (Section~\ref{SecA5}) suggests, reliable estimation of $p^*$ requires $1/\sqrt{N}$ to be of smaller order than $d^*$. However, for small $d^*$, linear predictors of orders $p^*-1$, $p^*$, and $p^*+1$ yield nearly identical final prediction errors, making exact identification of $p^*$ less relevant. 
    } 
\end{remark}

\subsection{Order selection by hypotheses testing}
\label{sec33} 

This section examines whether a linear predictor of prespecified order $p_0$ achieves a coefficient of determination of at least $\nu$. This can be assessed by testing
\begin{align}\label{det92}
    H_0 : p_0 < p^* \quad\text{ vs. }\quad H_1 : p_0 \geq p^*.  
\end{align}
By the definition of $p^*$ in \eqref{2.18}, the alternative is equivalent to both 
\[
   S_{p_0} \leq S_{p^*} \leq 1 - \nu  ,
    \quad\text{ and }\quad
    R^2_{p_0} \geq R^2_{p^*} \geq \nu.
\]
Thus, rejecting $H_0$ provides evidence at significance level $\alpha$ that a linear predictor of order $p_0$ achieves a coefficient of determination of at least $\nu$.

In the following, we develop a test for the hypotheses in~\eqref{det92}.
First, note that
\begin{align*}
    \mathbb{P}_{H_0} \big(\widehat{p} \leq p_0\big)  
    &=  1 - \mathbb{P}_{p^* > p_0} \Big(\bigcap^{p_0}_{p=1} \big\{\widehat{{T}}_p(\nu) \leq q_\alpha(W) \big\}\Big) ,    
\end{align*}
 where $\widehat{T}_p(\nu)$ is defined by
$$
\widehat T_p(\nu) \coloneqq\frac{1- \widehat{S}_p-\nu }{\widehat{V}_{S_p}} = \frac{S_p- \widehat S_p}{\widehat{V}_{S_p}} + \frac{1- S_p- \nu }{\widehat{V}_{S_p}}. 
$$
As $\widehat{T}_p(\nu)~ {\stackrel{\mathbb{P}}{\longrightarrow}} - \infty$, whenever $p < p^*,$ it follows that $ \mathbb{P}_{H_0} \big(\widehat{p} \leq p_0\big) \longrightarrow 0$.
This means that, under the decision rule {\it  ‘‘reject the null hypothesis in \eqref{det92} whenever $\widehat{p} \leq p_0$''} , the type~I error cannot be controlled. A valid test for \eqref{det92}, however, can be obtained by noting that these hypotheses are equivalent to  
\begin{equation}\label{2.17}
    H_0: S_{p_0} > 1-\nu \quad\text{ vs. }\quad H_1: S_{p_0} \leq 1-\nu,
\end{equation}
which were considered in Section~\ref{sec31}. The rejection rule
\begin{equation}\label{deq18}
    \widehat{S}_{p_0} \leq 1-\nu + q_{\alpha}(W)\,\widehat V_{S_{p_0}}
\end{equation}
defines an asymptotic level-$\alpha$ test
for \eqref{det92} and \eqref{2.17},
which is consistent whenever
$S_{p_0}<1-\nu$.

We conclude this section by testing whether a linear predictor of order $p_0$ does not attain a desired accuracy, through the hypotheses 
\begin{align}\label{det92'}
    H_0 : p_0 \geq p^* \quad\text{ vs. }\quad H_1 : p_0 < p^*.  
\end{align}
We consider the following test:
\begin{align}
\label{det82}
    \text{ reject } ~ H_0: p_0 \geq p^*, \text{  whenever ~} 
    \widehat p  > p_0.
\end{align}

\begin{theorem}\label{thm4}Under the assumptions of Theorem \ref{thm: Distr. free limit of M'}, the test \eqref{det82} has asymptotic level $\alpha$ and is consistent for the hypotheses in \eqref{det92'}.
\end{theorem}

\begin{figure}[h]
    \begin{minipage}{\textwidth}
        \centering
        \includegraphics[width=0.95\textwidth]{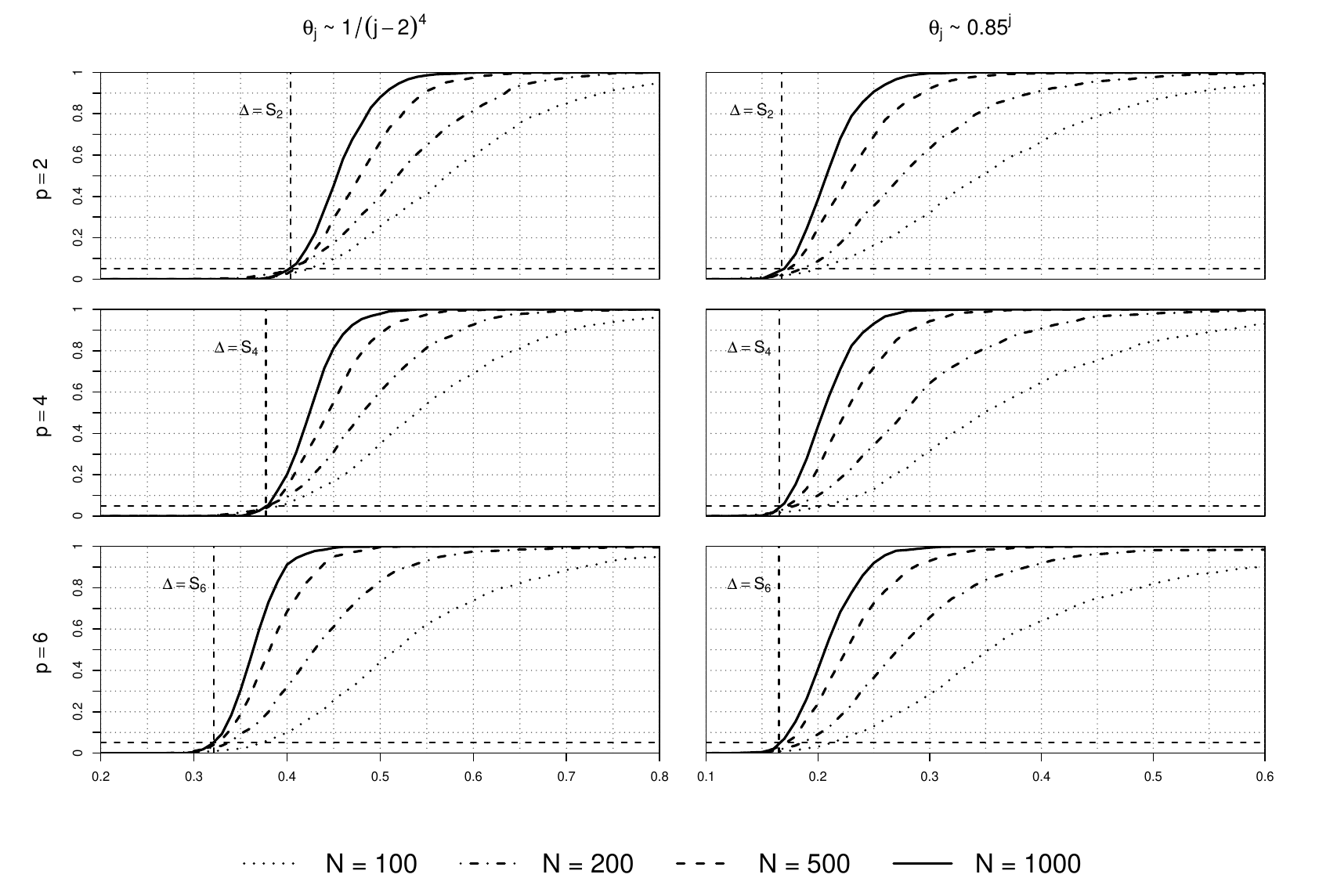}
    \end{minipage}   
\caption{\it Simulated rejection probabilities (y-axis) of the test  \eqref{det102} for the hypotheses \eqref{det103} for various values of the threshold $\Delta$ (x-axis). Vertical lines indicate the boundary of the hypotheses, where  $S_p=\Delta$, and horizontal lines mark the nominal level $\alpha = 5\%$. The data generating process is given by \eqref{eq:(X_k) lin process} with two choices of  coefficients given by \eqref{kue10}.}
\label{fig:uni} 
\end{figure}

\subsection{Finite sample properties}\label{sec34} 
We simulate the test \eqref{det102} and the estimation of the minimal lag $p^*$ for which the final prediction error $S_p$ in \eqref{det96} does not exceed $1-\nu$. Both experiments use $1000$ replications with i.i.d.\ innovations $\varepsilon_k\sim\mathcal{N}(0,1)$. The self-normalizer $\widehat{V}_{S_p}$ in \eqref{Vhat_S} is approximated by a Riemann sum with step size $1/20$,
starting at $1/20$ to avoid numerical instabilities.

\paragraph{Testing relevant hypotheses.} We begin with the test \eqref{det102} for the hypotheses \eqref{det103}. Assume that $(X_k)$ is the linear process in \eqref{eq:(X_k) lin process} with either polynomially or geometrically decaying coefficients:  
\begin{align}\label{kue10}
    \theta_j &= 
    \begin{cases}
        \;(j-2)^{-4}, & j>3,\\
        \;1,          & 0\leq j \leq 3,
    \end{cases}
    \qquad \text{or} \qquad 
    \theta_j = 
    \begin{cases}
        \;0.85^{j},  & j>3,\\ 
        \;2/3,       & 0\leq j \leq 3,
    \end{cases}
\end{align}
with $(S_2,S_4,S_6)\approx(0.404,0.377,0.322)$ for polynomial and $(0.167,0.163,0.161)$ for geometric decay. Figure~\ref{fig:uni} shows simulated rejection probabilities across thresholds $\Delta$ for the orders $p\in\{2,4,6\}$ and sample sizes $N\in\{100,200,500,1000\}$ at the nominal level $\alpha=5\%$. The results reflect Corollary~\ref{Cor:Asymptotic level test - distr. free}: as $N$ increases, rejection probabilities converge rapidly to 0 in the interior of the null hypothesis ($S_p>\Delta$), approach $\alpha$ at its boundary ($S_p=\Delta$), and 1 under the alternative ($S_p<\Delta$).

\begin{figure}[htbp]
    \begin{minipage}{\textwidth}
        \centering
        \includegraphics[width=0.9\textwidth]{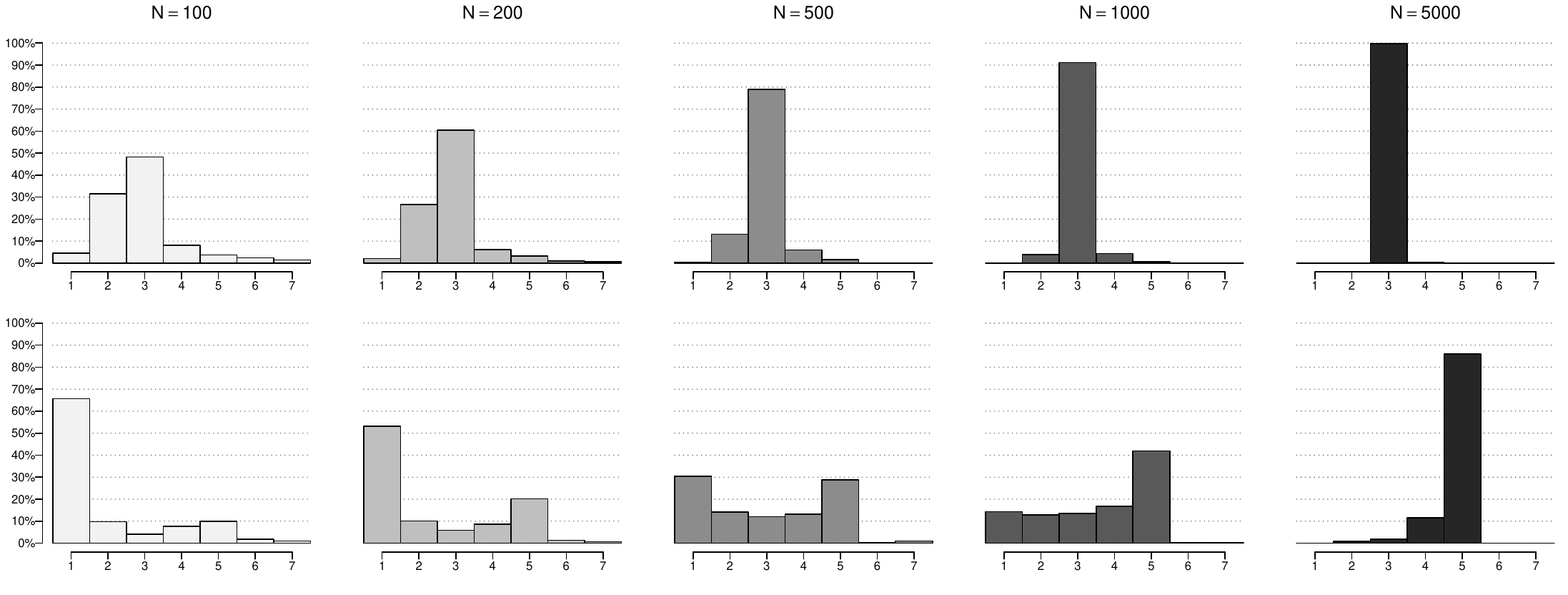}
    \end{minipage}
\caption{\it Histograms of the estimator $\widehat{p}$ for the lag $p^*$ defined in \eqref{2.18}, where  the nominal level is $\alpha = 10\%$. Upper panel: the process is given by the AR($5$) model defined in \eqref{kue13}. Since $\nu=0.6$ and thus $1-\nu = 0.4,$ the true value is given by $p^\ast = 3.$ Bottom panel: the process is given by the MA$(\infty)$ model defined in \eqref{kue10} with polynomially decaying coefficients. Here, $\nu=0.65$ and thus $1-\nu = 0.35,$ so the true value is $p^\ast = 5.$} 
\label{fig:hist1}
\end{figure}

\paragraph{Estimating the order for linear predictions.} 
We now evaluate the finite-sample performance of the estimator of $p^\ast=\min\{p\in\mathbb{N}\mid S_p\leq 1-\nu\}$.
First, we set $\nu=0.6$, and consider the AR(5) process
\begin{align}\label{kue13}
    X_k = -0.25X_{k-1} +0.1X_{k-2} + 0.4X_{k-3} - 0.25X_{k-4} +0.25X_{k-5} + \varepsilon_k, \qquad k \in \Z,
\end{align}
Table~\ref{tab:1} reports $S_p$ for $p=1,\ldots,7$, yielding $p^*=3$ and $S_{p^*}\approx0.366$. The upper part of Figure~\ref{fig:hist1} shows histograms of $\widehat p$ from \eqref{hx3} for various sample sizes $N$, with nominal level $\alpha=10\%$ to control overestimation of $p^*$. The estimator performs reasonably well, with accuracy improving as $N$ increases. The controlled probability of overestimating $p^*$ is clearly reflected in the results.

\begin{table}[H]
\centering
\caption{\it True $S_p$ values for the AR$(5)$ process \eqref{kue13} (top row) and the linear process with coefficients  \eqref{kue10} (bottom row).}
\label{tab:1} 
\begin{tabular}{r|c|c|c|c|c|c|c} \hline $p$ & 1 & 2 & 3 & 4 & 5 & 6 & 7 \\ 
    \hline 
    AR$(5)$ & 0.679 & 0.613 & 0.366 & 0.325 & 0.305 & 0.305 & 0.305\\
    \hline
    MA$(\infty)$ & 0.415 & 0.404 & 0.393 & 0.377 & 0.324 & 0.322 & 0.318\\
\hline 
\end{tabular}
\end{table}

\medskip

\noindent To illustrate the conceptual difference between our approach and standard model selection criteria, consider the AIC, which balances goodness of fit and model complexity. For $N=500$ observations from the AR$(5)$ process \eqref{kue13}, the AIC selects $p=5$ in $70\%$, $p=6$ in $10\%$, and $p\geq 7$ in $20\%$ of cases among AR models with $p\leq 9$, all yielding $S_p=0.305$. By contrast, our pivotal method targets the smallest lag order satisfying $S_{p^*}\leq 0.4$, namely $p^*=3$. Thus, rather than selecting a specific model, it determines the predictor order needed to achieve a prespecified prediction accuracy.

Next, we set $\nu=0.65$ and consider the linear process with polynomially decaying coefficients in \eqref{kue10}. The true $S_p$ values (Table~\ref{tab:1}, bottom row) give $p^*=5$ for $1-\nu=0.35$. The histograms of $\widehat p$ (Figure~\ref{fig:hist1}, bottom panel)
show lower accuracy than for the AR$(5)$ process, reflecting
the smaller value of
$d^*=\max\{S_{p^*-1}-S_{p^*},S_{p^*}-S_{p^*+1}\}$:
$0.053$ versus $0.247$ for AR$(5)$.
Nevertheless, overestimation of $p^*$ remains rare, as intended
by our procedure. By contrast, AIC consistently selects the
maximum candidate order $p=9$ among AR models with $p\leq 9$,
reflecting model misspecification.

\section{Relative improvement and partial autocorrelations}\label{sec4}

This section considers the measure $Q_p$ in \eqref{det91},
which compares the final prediction errors from linear
predictors of orders $p$ and $p-1$. We define the
corresponding sequential estimator by
\begin{align*}
    \widehat Q_p(\lambda)\coloneqq
    \begin{cases}
        \,\displaystyle
        \frac{\widehat M_p(\lambda)}
        {\widehat M_{p-1}(\lambda)},
        & \lambda\in(0,1],\;
        \widehat M_{p-1}(\lambda)>0,\\[6pt]
        \,0,
        & \lambda\in(0,1],\;
        \widehat M_{p-1}(\lambda)=0,\\
        \,1, & \lambda=0,
    \end{cases}
\end{align*}
with $\widehat Q_p\coloneqq\widehat Q_p(1)$
denoting the full-sample estimator of $Q_p$.
For completeness, we set $\widehat{Q}_p(0) \coloneqq 1$ and introduce the statistic
\begin{align*}
    \widehat{V}_{Q_p} \coloneqq \int^1_0 \lambda\big|\widehat{Q}_p(\lambda) -  \widehat{Q}_p\big|\,\mathrm{d}\lambda.
\end{align*}

\begin{theorem}\label{thm:piv_Q} 
Suppose $\Sigma_p$ in \eqref{AN for vector of autocov} is non-singular and $\kappa_p \neq 0$. Then, with $W$ as in \eqref{def W},
\begin{align*}
    \frac{\widehat{Q}_p-Q_p}{\widehat{V}_{Q_p}}\, \,\stackrel{d}{\longrightarrow} \,W.
\end{align*}
\end{theorem}

An asymptotic $(1-\alpha)$ confidence interval for $Q_p$ is thus given by 
\begin{align*}
    \Big[\,\widehat Q_p - q_{1-\alpha/2} \big( W\big) {\widehat{V}_{Q_p}}~,~ \widehat Q_p + q_{1-\alpha/2}\big(W \big) {\widehat{V}_{Q_p}}\Big].
\end{align*}
Similarly, a pivotal,  consistent and asymptotic level $\alpha$-test for the hypotheses
\begin{align*}
    H_0 : Q_p \geq \Delta \quad\text{ vs. }\quad H_1 : Q_p < \Delta,
\end{align*}
is obtained by rejecting the null hypothesis, whenever $\widehat{Q}_p \leq \Delta +  q_{\alpha}(W)\widehat{V}_{Q_p}.$
 
\subsection{Partial autocorrelations}\label{sec42}

As noted in the introduction, our results provide new tools for inference on partial autocorrelations, playing a central role in selecting the order of autoregressive (AR) models \citep[see][]{Durbin1960}, since $\kappa_h=0$ for all $h>p$ in an AR($p$) process. It follows from \citet[][Section~5.2]{BrockwellDavis1991} that
\begin{align}\label{relationship_kappa_Q}
    1-\kappa_h^2 = Q_h = \frac{M_h}{M_{h-1}}, \qquad h \in \N .
\end{align}
The decision rule $\widehat{Q}_{p_0} \leq  1- \nu  + q_{\alpha}  (W) \widehat V_{Q_{p_0}}$ provides an asymptotic level-$\alpha$ test of whether $\kappa_{p_0}^2$ exceeds $\nu.$ Moreover, a pivotal confidence interval for the squared partial autocorrelation $\kappa_p^2>0$ is given by
\begin{align*}
    \Big[\widehat \kappa_p^2 - q_{1-\alpha/2}(W)\widehat{V}_{Q_p} ~,~ \widehat \kappa_p^2 + q_{1-\alpha/2}(W)\widehat{V}_{Q_p}\Big],
\end{align*}
where $\widehat \kappa_p^2 = 1 - \widehat Q_p$. Using the representation
\[
    \kappa_p = \boldsymbol{e}^\top_p G^{-1}_{p-1} (\gamma_1, \dots, \gamma_p)^\top,
\]
with $G_{p-1}=(\gamma_{j-i})_{i,j=0}^{p-1},$ from the Durbin–Levinson algorithm \citep[][Eq.\ (3.4.2)]{BrockwellDavis1991} and $\boldsymbol{e}_p$ the $p$th unit vector in $\R^p$, we can further construct a pivotal estimator of $\kappa_p$. 
Specifically, we define the sequential estimator
\begin{align}\label{kue5}
    \widehat\kappa_p(\lambda)\coloneqq
    \begin{cases}
        \,\displaystyle
        \boldsymbol e_p^\top
        \widehat G_{p-1}(\lambda)^{-1}
        \big(\widehat\gamma_1(\lambda),\dots,
        \widehat\gamma_p(\lambda)\big)^\top,
        & \lambda\in(0,1],\;
        \widehat G_p(\lambda)\succ0,\\[6pt]
        \,0, & \text{otherwise},
    \end{cases}
\end{align}
where
$\widehat G_k(\lambda)
=(\widehat\gamma_{|j-i|}(\lambda))_{i,j=0}^k$,
$k=p-1,p$, and $\boldsymbol e_p$ denotes
the $p$th unit vector in $\mathbb R^p$.
We further write
$\widehat\kappa_p\coloneqq\widehat\kappa_p(1)$
for the full-sample estimator.
and define
\begin{align}\label{Vhat_K}
    \widehat{V}_{\kappa_p} \coloneqq \int_0^1 \lambda \big|\widehat{\kappa}_p(\lambda) - \widehat{\kappa}_p\big|\,\mathrm{d}\lambda\,.   
\end{align}

\begin{theorem}\label{thm:PACs} Let $\Sigma_p$ in \eqref{AN for vector of autocov} be non-singular. Then, with $W$ as in \eqref{def W},
\begin{align*}
    \frac{\widehat{\kappa}_p-\kappa_p}{\widehat{V}_{\kappa_p}}\, \,\stackrel{d}{\longrightarrow} \,W.
\end{align*}
\end{theorem}

As an immediate consequence of Theorem~\ref{thm:PACs} we obtain the pivotal confidence interval for $\kappa_p$,
\begin{align}\label{kue16}
    \Big[\widehat\kappa_p - q_{1-\alpha/2}(W)\widehat V_{\kappa_p} ~,~ \widehat\kappa_p + q_{1-\alpha/2}(W)\widehat V_{\kappa_p}\Big].
\end{align}
For comparison, applying the delta method to \eqref{AN for vector of autocov} yields
\begin{align*}
    \sqrt{N}(\widehat\kappa_p-\kappa_p) \stackrel{d}{\longrightarrow} \mathcal{N}(0,\theta_p).
\end{align*}
By \eqref{Def of Sigma_ij}, the asymptotic variance $\theta_p$ has a Bartlett-type structure involving all autocovariances $\gamma_k$, making estimation difficult \citep[cf.][for a Bartlett-type formula with Gaussian innovations]{Stoica1989}. For AR($p$) processes with Gaussian innovations, $\theta_p=1-\kappa_p^2$ \citep{Barndorff-NielsenSchou1973}, but the asymptotic variances of $\widehat\kappa_1,\ldots,\widehat\kappa_{p-1}$ remain difficult to estimate. By contrast, self-normalization yields pivotal, asymptotically valid confidence intervals for partial autocorrelations of any order in a linear process.

\subsection{Finite sample properties}\label{sec43}

We investigate the finite-sample performance of the pivotal confidence intervals \eqref{kue16}, comparing our method (PIV) with that of \cite{Barndorff-NielsenSchou1973} (BNS). The latter estimates AR coefficients by maximum likelihood under a specified model order and transforms them into partial autocorrelations. As before, we assume i.i.d.\ innovations $\varepsilon_k\sim\mathcal{N}(0,1)$.

Confidence intervals for $\kappa_2$ and $\kappa_4$ are considered under two scenarios. (i) BNS assumes the correct AR order ($p=2$ or $p=4$; left panel of Table~\ref{tab:comparison_covering_length}), where the AR(2) and AR(4) models are given by 
\begin{align*}
    X_k = -0.2 X_{k-1} - 0.3 X_{k-2} + \varepsilon_k, \qquad
    X_k = -0.2 X_{k-1} - 0.3 X_{k-2} + 0.3 X_{k-3} + 0.2 X_{k-4} + \varepsilon_k,
\end{align*}
with $\kappa_2 = -0.3$ and $\kappa_4 = 0.2$, respectively. (ii) BNS incorrectly fits AR(2) and AR(4) models (right panel of Table \ref{tab:comparison_covering_length}), while the data are generated from the AR(6) process
\begin{align*}
    X_k = -0.2 X_{k-1} - 0.3 X_{k-2} + 0.3 X_{k-3} + 0.2 X_{k-4} + 0.1 X_{k-5} + 0.1 X_{k-6} + \varepsilon_k,
\end{align*}
in which $\kappa_2\approx -0.377$ and $\kappa_4\approx 0.157$. We consider the sample sizes $N \in \{100,200,500,1000\}$, nominal level $\alpha=10\%$, and $\widehat{V}_{\kappa_p}$ in~\eqref{Vhat_K} is computed by a Riemann sum with step size $1/20$, starting at $1/20$. 

The finite-sample results in Table~\ref{tab:comparison_covering_length} reflect Theorem~\ref{thm:PACs}. In scenario (i), both methods achieve coverage near the nominal level $1-\alpha=0.9$, with interval lengths decreasing as $N$ grows. PIV intervals are slightly wider, reflecting a moderate efficiency loss in exchange for greater robustness, as noted by \cite{shao2015}. This advantage becomes evident in scenario (ii), where the BNS assumptions are violated. While PIV maintains similar performance, BNS coverage drops to around 80\% for $\kappa_2$ and slightly below 90\% for $\kappa_4$, despite shrinking interval lengths. Thus, PIV performs comparably under correct specification and substantially more reliably under misspecification.

\begin{table}[ht]
\centering
\caption{Comparison of our pivotal (PIV) confidence intervals \eqref{kue16} with those of \cite{Barndorff-NielsenSchou1973} (BNS) for $\kappa_2,\kappa_4$ at confidence level $1-\alpha=0.9$, showing empirical coverage and interval length. The left panel shows scenario\,(i) with correctly specified AR orders ($p\in\{2,4\}$), and the right scenario\,(ii), fitting AR(2) or AR(4) models to an AR(6) process.
}
\label{tab:comparison_covering_length} 
\resizebox{0.925\textwidth}{!}{%
\begin{tabular}{c@{\hskip 0.5in}c}
\textbf{(i) Correct order by BNS} & \textbf{(ii) Incorrect order by BNS} \\

\begin{tabular}{cr|cc|cc}
\hline
\multirow{2}{*}{$ p $} & \multirow{2}{*}{$ N $} 
& \multicolumn{2}{c|}{BNS} 
& \multicolumn{2}{c}{PIV} \\
& & {\footnotesize Coverage} & {\footnotesize Length} & {\footnotesize Coverage} & {\footnotesize Length} \\
\hline
\multirow{4}{*}{2} 
& 100        & 0.902 & 0.314 & 0.893 & 0.386 \\
& 200        & 0.911 & 0.222 & 0.903 & 0.281 \\
& 500        & 0.908 & 0.140 & 0.898 & 0.181 \\
& 1000 & 0.898 & 0.099 & 0.903 & 0.129 \\
\hline
\multirow{4}{*}{4} 
& 100         & 0.880 & 0.329 & 0.925 & 0.515 \\
& 200         & 0.909 & 0.230 & 0.908 & 0.331 \\
& 500         & 0.896 & 0.145 & 0.909 & 0.189 \\
& 1000  & 0.888 & 0.102 & 0.898 & 0.131 \\
\hline
\end{tabular}
& \hspace{-0.5cm}
\begin{tabular}{cr|cc|cc}
\hline
\multirow{2}{*}{$ p $} & \multirow{2}{*}{$ N $} 
& \multicolumn{2}{c|}{BNS} 
& \multicolumn{2}{c}{PIV} \\
& & {\footnotesize Coverage} & {\footnotesize Length} & {\footnotesize Coverage} & {\footnotesize Length} \\
\hline
\multirow{4}{*}{2} 
& 100        & 0.788 & 0.300 & 0.872 & 0.451 \\
& 200        & 0.813 & 0.215 & 0.897 & 0.335 \\
& 500        & 0.786 & 0.136 & 0.893 & 0.221 \\
& 1000 & 0.790 & 0.096 & 0.901 & 0.158 \\
\hline
\multirow{4}{*}{4} 
& 100        & 0.865 & 0.332 & 0.913 & 0.495 \\
& 200        & 0.892 & 0.232 & 0.908 & 0.317 \\
& 500        & 0.878 & 0.146 & 0.913 & 0.197 \\
& 1000 & 0.876 & 0.103 & 0.905 & 0.134 \\
\hline
\end{tabular}
\end{tabular}
}
\end{table}

\section{Multivariate setting}\label{sec5}

This section illustrates our approach for multivariate stationary processes $(X_k)_{k\in\Z}\subset\R^d$ of the form 
\begin{align}\label{eq:(X_k) lin process mult}
    X_k = \big ( X_k^{(1)}, \ldots,  X_k^{(d)} \big )^\top = \sum^\infty_{j=0}\,\Theta_{\!j}\varepsilon_{k-j},\qquad k\in\Z,  
\end{align}
where $d\in\N$, $\Theta_{\!j}\in\R^{d\times d}$ are matrices with $\sum_{j=1}^\infty j\|\Theta_j\|<\infty$, and $(\varepsilon_k)_{k\in\Z}$ is a sequence of i.i.d.\ $d$-dimensional random variables with $\Erw(\varepsilon_0)=0$ and $\Erw(\varepsilon_0\varepsilon_0^\top)=I_d$ (the $d\times d$ identity matrix). We further assume that the components of the innovations $\varepsilon_0$ have finite fourth moments and denote by $\Gamma_{\!h}=\Gamma^\top_{\!-h}=\Erw(X_0X_h^\top)$ the corresponding autocovariance matrices. 

To define an analogue of $M_p$ in \eqref{Mp}, let $\|\cdot\|_2$ denote the Euclidean norm on $\R^d$ and consider a linear predictor of the form $\sum_{i=1}^p \Xi_i X_{k-i}$ with coefficient matrices $\Xi_i \in \R^{d\times d}.$ A general result on linear approximation in Hilbert spaces \citep[see][p.\ 16]{Achieser1956} then yields
\begin{align}
    \mathscr{M}_p 
    \coloneqq 
    \underset{\Xi_1, \dots, \Xi_p\,\in\,\R^{d\times d}}{\min}\Erw\!\Big\|X_p - \sum^p_{i=1}\,\Xi_iX_{p-i}\Big\|^2_2
    \,=\, 
    \frac{1}{\det(\boldsymbol{\mathscr{G}}_{\!p-1})}\sum^d_{j=1}\,\det(\boldsymbol{\mathscr{G}}_{\!p-1,j}),  \label{det200}
\end{align}
where $\boldsymbol{\mathscr{G}}_{\!p-1}\coloneqq(\Gamma_{\!j-i})^{p-1}_{i,j=0}\in\R^{dp\times dp}$ is the block autocovariance matrix, assumed to be non-singular throughout this section, $\boldsymbol{\mathscr{G}}_{\!p-1,j}\in\R^{(dp+1)\times(dp+1)}$ are matrices defined by
\[
    \boldsymbol{\mathscr{G}}_{\!p-1,j} \coloneqq 
\left(
\begin{NiceArray}{c|cccc}
\boldsymbol{e}^\top_j \Gamma_{\!0} \boldsymbol{e}_j
  & \boldsymbol{e}^\top_j \Gamma^\top_{\!1} 
  & \cdots 
  & \boldsymbol{e}^\top_j \Gamma^\top_{\!p} \\
\hline
\Gamma_{\!1} \boldsymbol{e}_j 
  & \Block{3-3}{\boldsymbol{\mathscr{G}}_{\!p-1}} \\
\vdots \\
\Gamma_{\!p} \boldsymbol{e}_j \\
\end{NiceArray}
\right)
,\qquad j = 1,\dots,d,
\]
and $\boldsymbol{e}_j$ is the $j$th unit vector in $\R^d\!.$ Also, we define $\mathscr{M}_0 = \Erw\!\| X_0\|_2^2 = \operatorname{tr}(\Gamma_{\!0})$.

For the sake of brevity we restrict ourselves to  the normalized measure
\begin{align}\label{det84}
    \mathscr{S}_p \coloneqq  \frac{\mathscr{M}_p}{\mathscr{M}_0}=\frac{\mathscr{M}_p}{\tr(\Gamma_{\!0})},      
\end{align}
which defines a multivariate analogue of $S_p$ from Section~\ref{sec3}. Results for $\mathscr{M}_p/\mathscr{M}_{p-1}$ from Section~\ref{sec4} follow analogously.

As in Section~\ref{sec2}, we define the sequential
estimator of $\mathscr M_p$ by
\begin{align}\label{det201}
    \widehat{\mathscr M}_p(\lambda)
    \coloneqq
    \begin{cases}
    \,\displaystyle
    \frac{
    \sum_{j=1}^d
    \det(\widehat{\boldsymbol{\mathscr G}}_{\!p-1,j}(\lambda))
    }{
    \det(\widehat{\boldsymbol{\mathscr G}}_{\!p-1}(\lambda))
    },
    & \widehat{\boldsymbol{\mathscr G}}_{\!p}(\lambda)
    \succ0,\\[8pt]
    \,0, & \text{otherwise},
    \end{cases}
    \qquad \lambda\in[0,1],
\end{align}
where
$\widehat{\mathscr M}_p
\coloneqq\widehat{\mathscr M}_p(1)$.
Here,
$\widehat{\boldsymbol{\mathscr G}}_{\!p}(\lambda)
=(\widehat\Gamma_{\!j-i}(\lambda))_{i,j=0}^p$,
and the matrices
$\widehat{\boldsymbol{\mathscr G}}_{\!p-1}(\lambda)$
and
$\widehat{\boldsymbol{\mathscr G}}_{\!p-1,j}(\lambda)$
are obtained from
$\boldsymbol{\mathscr G}_{\!p-1}$ and
$\boldsymbol{\mathscr G}_{\!p-1,j}$,
respectively, by replacing the autocovariance
matrices $\Gamma_{\!h}$ with the sequential counterparts
\[
\widehat\Gamma_{\!h}(\lambda)
\coloneqq
\frac{1}{N}
\sum_{i=1}^{\lfloor\lambda(N-h)\rfloor}
X_iX_{i+h}^\top,
\qquad
0\leq h<N,\quad\lambda\in[0,1],
\]
and
$\widehat\Gamma_{\!h}(\lambda)
\coloneqq\widehat\Gamma_{\!-h}^\top(\lambda)$
for $-N<h<0$. For $p=0$, we set
$\widehat{\mathscr M}_0(\lambda)
\coloneqq\operatorname{tr}(\widehat\Gamma_0(\lambda))$. Further, we introduce the self-normalizer
\begin{align}\label{kue14}
\widehat{\mathscr V}_{\mathscr S_p}
\coloneqq
\int_0^1
\lambda\big|
\widehat{\mathscr S}_p(\lambda)
-\widehat{\mathscr S}_p
\big|\,\mathrm d\lambda,
\end{align}
where the sequential estimator is defined by
\[
    \widehat{\mathscr S}_p(\lambda)
    \coloneqq
    \begin{cases}
        \,\displaystyle
        \frac{\widehat{\mathscr M}_p(\lambda)}
        {\widehat{\mathscr M}_0(\lambda)},
        & \lambda\in(0,1],\;
        \widehat{\mathscr M}_0(\lambda)>0,\\[6pt]
        \,0,
        & \lambda\in(0,1],\;
        \widehat{\mathscr M}_0(\lambda)=0,\\
        \,1, & \lambda=0,
    \end{cases}
\]
and
$\widehat{\mathscr S}_p
\coloneqq\widehat{\mathscr S}_p(1)$
denotes the full-sample estimator of
$\mathscr S_p$ in \eqref{det84}.

Let $\mathrm{vech}(\cdot)$ and $\mathrm{vec}(\cdot)$ stack the lower triangular and all matrix entries columnwise, respectively. Since $\Gamma_{\!0}$ is symmetric but $\Gamma_{\!h}$ for $h\neq 0$ need not be, we define
\[
    \boldsymbol{\nu}_p \coloneqq
    \big(\mathrm{vech}^\top(\Gamma_{\!0}),
    \mathrm{vec}^\top(\Gamma_{\!1}),\dots,
    \mathrm{vec}^\top(\Gamma_{\!p})\big)^\top
    \in\R^{k_{p,d}},
\]
where $k_{p,d}\coloneqq d(d+1)/2+pd^2$, and define $\widehat{\boldsymbol{\nu}}_p$ analogously, with $\widehat\Gamma_{\!h}\coloneqq\widehat\Gamma_{\!h}(1)$ for $h=0,\dots,p$. Then,
\begin{align}\label{det203}
\sqrt{N}\big(\widehat{\boldsymbol{\nu}}_p-\boldsymbol{\nu}_p\big)
\stackrel{d}{\longrightarrow}
\mathcal{N}(0,\boldsymbol{\Sigma}_p),
\end{align}
where $\boldsymbol{\Sigma}_p\in\R^{k_{p,d}\times k_{p,d}}$. Since $\mathscr{M}_p$ in \eqref{det200} depends on $\boldsymbol{\nu}_p$,  $\mathscr{S}_p$ in \eqref{det84} can, for an appropriate function $g:\R^{k_{p,d}}\to\R,$ be expressed as
\begin{align}\label{det202}
    \mathscr{S}_p = g(\boldsymbol{\nu}_p).
\end{align}

\begin{theorem}\label{piv_M_cor mult}
Suppose $\boldsymbol{\Sigma}_p$ in \eqref{det203} is non-singular and $\nabla g(\boldsymbol{\nu}_p)\neq 0$. Then, with $W$ defined in \eqref{def W},
\[
    \frac{\widehat{\mathscr{S}}_p-\mathscr{S}_p}
    {\widehat{\mathscr{V}}_{\mathscr{S}_p}}
    \stackrel{d}{\longrightarrow} W.
\]
\end{theorem}

Several statistical applications can be derived in a similar manner as described in Sections \ref{sec31}--\ref{sec33}. Exemplarily, we propose a test for the hypotheses
\begin{align}\label{det85}
    H_0 : \mathscr{S}_p \geq \Delta \quad\text{ vs. }\quad H_1 : \mathscr{S}_p < \Delta , 
\end{align}
which rejects the null hypothesis,  whenever
\begin{align}\label{det86}
    \widehat{\mathscr{S}}_p \leq \Delta +  q_{\alpha}(W)\widehat{\mathscr{V}}_{\mathscr{S}_p}.
\end{align}
As in Section~\ref{sec31}, this decision rule yields a pivotal, consistent test with asymptotic level $\alpha$.

\begin{figure}[h]
\begin{minipage}{\textwidth}
\centering
\includegraphics[width=0.95\textwidth]{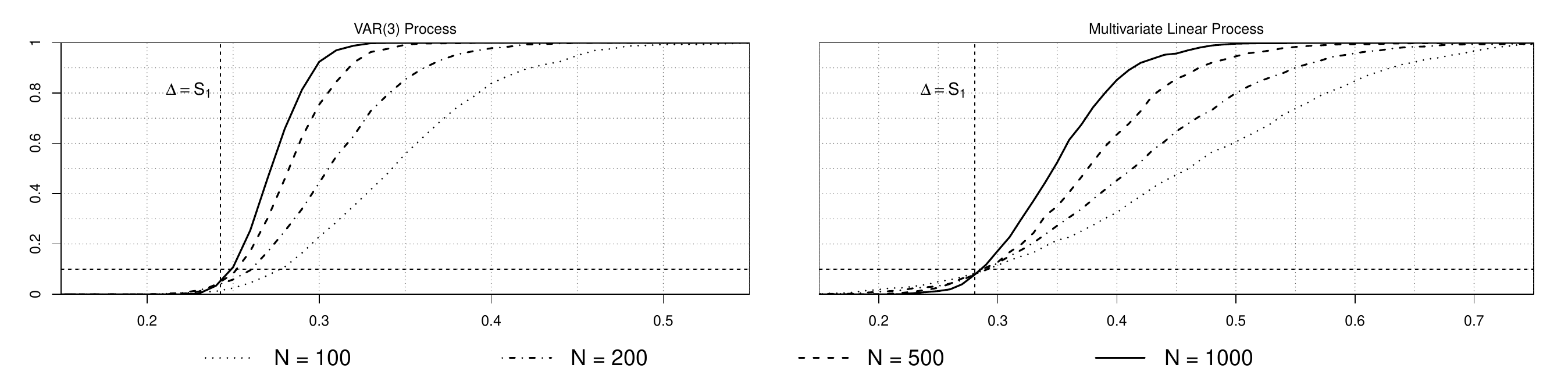}
\end{minipage}   
\caption{\it Empirical rejection probabilities ($y$-axis) against $\Delta$ ($x$-axis) for the VAR$(3)$ process (left) and the multivariate linear process (right). Vertical lines mark the true $\mathscr{S}_1$ values, and the horizontal lines the nominal level $\alpha=10\%$.}
\label{fig:Approx_multi}
\end{figure}

We illustrate its finite-sample performance using two $5$-dimensional stationary processes. The first is the VAR$(3)$ process
\begin{align}\label{det300}
    X_k=\Phi_1X_{k-1}+\Phi_2X_{k-2}
    +\Phi_3X_{k-3}+\varepsilon_k,\qquad k\in\Z,
\end{align}
where the innovations $\varepsilon_k$ are i.i.d.\
$\mathcal N_5(0,I_5)$ and the coefficient matrices
are given by
$$
\Phi_1 \coloneqq \scalebox{0.95}{$0.16\cdot\begin{pmatrix}
7 & 2 & 1 & 0 & 0   \\
2 & 5 & 2 & 1 & 0   \\
1 & 2 & 5 & 2 & 1 \\
0 & 1 & 2 & 5 & 2 \\
0 & 0 & 1 & 2 & 5 \\
\end{pmatrix}
$
}\!, \qquad
\Phi_2 \coloneqq \scalebox{0.95}{$-0.1\cdot
\begin{pmatrix}
3 & 2 & 0 & 0 & 0 \\
2 & 3 & 2 & 0 & 0 \\
0 & 2 & 3 & 2 & 0 \\
0 & 0 & 2 & 3 & 2 \\
0 & 0 & 0 & 2 & 3 \\
\end{pmatrix}
$
}\!, \qquad
\Phi_3 \coloneqq \scalebox{0.95}{$-0.05\cdot
\begin{pmatrix}
    2 & 1 & 0 & 0 & 0 \\
    1 & 1 & 1 & 0 & 0 \\
    0 & 1 & 1 & 1 & 0 \\
    0 & 0 & 1 & 1 & 1 \\
    0 & 0 & 0 & 1 & 1 \\
\end{pmatrix}
.$
}
$$
Secondly, consider the linear process in~\eqref{eq:(X_k) lin process mult} with coefficient matrices
\[
    \Theta_j \coloneqq \left(\tfrac{3}{5}\Phi_1\right)^{j}\!, \qquad j \geq 0,
\]
where $\Phi_1$ is defined above. Figure~\ref{fig:Approx_multi} shows the simulated rejection probabilities for the test \eqref{det86} across thresholds $\Delta$ at level $\alpha=10\%$. The self-normalizer in \eqref{kue14} is approximated by a Riemann sum with step size $1/20$, starting at $1/20$, and autocovariances are computed using \texttt{VARMAcov()} from the \texttt{R} package \texttt{MTS} \citep{TsayWoodLachmannMTS2022}. The results mirror the univariate case in Figure~\ref{fig:uni}.

\paragraph{Data Availability Statement} The code used to reproduce all tables and figures reported in this article is available at \url{https://github.com/SebastianKuehnert/PivInf_LinPred_StatProc}.

\paragraph{Acknowledgements} This work was supported by SFB 1597 (Project-ID 499552394) and TRR 391 (Project-ID 520388526) funded by the Deutsche Forschungsgemeinschaft (DFG, German Research Foundation). The authors thank the anonymous referees for their valuable comments on earlier versions of this article, as well as their colleagues Patrick Bastian, Marius Kroll, and Zihao Yuan for helpful discussions.

\paragraph{Disclosure of AI Tool Use} During the preparation of this manuscript, ChatGPT (OpenAI) was used to assist with brainstorming, checking mathematical arguments, language editing, LaTeX formatting, and code development and optimization. The research questions, methodological framework, and principal scientific contributions were developed by the authors. All mathematical arguments, computational implementations, and results were implemented, and verified by the authors, incorporating AI-generated suggestions. All such suggestions were critically evaluated and revised where necessary. The authors take full responsibility for the work. No sensitive data were shared with the AI tool.

{\footnotesize
    \bibliographystyle{chicago}
    \bibliography{References_AR_MA_testing}
}

\appendix

\section{Appendix}

Before we prove the main results, we establish the following.

\begin{lemma}\label{LemA1} The gradient of the measure $M_p$ in \eqref{Def M_p}, with $p\in\mathbb{N}_0,$ evaluated at $c\boldsymbol{\gamma}_p$ satisfies
\begin{align}\label{grad_identity}
    \nabla M_{p,c\boldsymbol{\gamma}_p}
    =
    \nabla M_{p,\boldsymbol{\gamma}_p}, \qquad c \neq 0,
\end{align}
where $\boldsymbol{\gamma}_p = (\gamma_0, \dots, \gamma_p)^\top$ is the vector of autocovariances of the process in \eqref{eq:(X_k) lin process}. Moreover, $\nabla M_{p,\boldsymbol{\gamma}_p}\neq0$. 
\end{lemma}

\begin{proof} First, consider $p=0.$ Then, $M_p=M_0=\gamma_0,$ and the gradient of $M_0$ evaluated at any $x\in \mathbb{R}$ is $\nabla M_{0,x} = 1.$ Thus, for $p=0,$ it holds both $\eqref{grad_identity}$ and $\nabla M_{p,\boldsymbol{\gamma}_p} \neq 0.$

Next, let $p\in \mathbb{N}.$ As covariance matrices, the Toeplitz matrices $G_p=(\gamma_{j-i})_{i,j=0}^p$ are non-negative definite and, by our non-singularity assumption, positive definite for every \(p\). Further, we can write
\[
    G_p=
    \begin{pmatrix}
    G_{p-1} & \boldsymbol r_p\\
    \boldsymbol r_p^\top & \gamma_0
    \end{pmatrix},
    \qquad
    \boldsymbol r_p
    =(\gamma_p,\ldots,\gamma_1)^\top.
\]
The Schur complement formula and positive definiteness then yield
\begin{align}\label{Schur}
    M_p = \frac{\det(G_p)}{\det(G_{p-1})} = \gamma_0 - \boldsymbol r_p^\top G_{p-1}^{-1}\boldsymbol r_p > 0,
\end{align}
Since $G_p(c\boldsymbol\gamma_p)=cG_p(\boldsymbol\gamma_p)$ and $M_p(c\boldsymbol\gamma_p)=cM_p(\boldsymbol\gamma_p)$ for $c\neq0$, and since $\partial_{\gamma_\ell}G_p$ is independent of $\boldsymbol\gamma_p$, combining the quotient rule with Jacobi's formula gives for any $0 \leq \ell \leq p,$
\begin{align*}
    \partial_{\gamma_\ell}M_p(c\boldsymbol\gamma_p)
    &=
    M_p\Big(
    \operatorname{tr}\big(G_p^{-1}\partial_{\gamma_\ell}G_p\big)
    -
    \operatorname{tr}\big(G_{p-1}^{-1}\partial_{\gamma_\ell}G_{p-1}\big)
    \Big)\\
    &=
    cM_p\Big(
    \operatorname{tr}\big(c^{-1}G_p^{-1}\partial_{\gamma_\ell}G_p\big)
    -
    \operatorname{tr}\big(c^{-1}G_{p-1}^{-1}\partial_{\gamma_\ell}G_{p-1}\big)
    \Big)\\
    &=
    \partial_{\gamma_\ell}M_p(\boldsymbol\gamma_p).
\end{align*}
Hence, \eqref{grad_identity} is also satisfied for $p\in\mathbb{N}.$ Further, by the block inverse formula,
\[
    G_p^{-1}
    =
    \begin{pmatrix}
    G_{p-1}^{-1}
    +\dfrac{\boldsymbol u_p\boldsymbol u_p^\top}{M_p}
    &
    -\dfrac{\boldsymbol u_p}{M_p}
    \\[8pt]
    -\dfrac{\boldsymbol u_p^\top}{M_p}
    &
    \dfrac{1}{M_p}
    \end{pmatrix},
\]
where $\boldsymbol u_p=G_{p-1}^{-1}\boldsymbol r_p$. Since $\partial_{\gamma_0}G_p=I_{p+1}$ and $\partial_{\gamma_0}G_{p-1}=I_p$, where $I_k$ denotes the $k$-dimensional identity matrix, we obtain
\begin{align}
    \partial_{\gamma_0}M_p(\boldsymbol\gamma_p)
    \,=\,
    M_p\big(
    \operatorname{tr}(G_p^{-1})
    -
    \operatorname{tr}(G_{p-1}^{-1})
    \big)
    \,=\,
    1+\|G_{p-1}^{-1}\boldsymbol r_p\|^2.
    \label{part_der_Mp_gamma0}
\end{align}
Thus, we also have $\nabla M_{p,\boldsymbol\gamma_p}\neq0$ for $p\in\mathbb{N},$ which completes the proof.
\end{proof}

\subsection{Proof of Theorem \ref{M0, ..., Mp cvgc}}\label{SecA1}
We begin with introducing some notation. Recall that $\ell^{\infty,p+1}([0,1])$ is the space of bounded functions $\boldsymbol{f}_{\!p} = (f_0, \dots, f_p)^\top\colon[0,1]\rightarrow \R^{p+1}$ with norm $\|\boldsymbol{f}_{\!p}\|_\infty = \sup_{\lambda\in[0,1]}\sup_{0\leq i \leq p}|f_i(\lambda)|,$ and $\ell([0,1]) = \ell^{\infty,1}([0,1]),$ and $\ell^{\infty, p+1}_0([0,1]) = \{\boldsymbol{f}_{\!p} \in \ell^{\infty, p+1}([0,1]) \,|\, \boldsymbol{f}_{\!p}(0) = \boldsymbol{0}\}.$ Further, $\rightsquigarrow$ denotes weak convergence in such spaces. Moreover, we write
\begin{align}\label{Def hatg_p and g_p}
    \widehat{\boldsymbol{g}}_p&\coloneqq \{\widehat{\boldsymbol{\gamma}}_{\!p}(\lambda)\}_{\lambda\in[0,1]}\in\ell^{\infty,p+1}_0([0,1])\quad\;\text{and}\quad\;\boldsymbol{g}_p\coloneqq \{\lambda\boldsymbol{\gamma}_p\}_{\lambda\in[0,1]}\in\ell^{\infty,p+1}_0([0,1]), 
\end{align}
with $\boldsymbol{\gamma}_p = (\gamma_0, \ldots , \gamma_p)^\top\!,$ and $\widehat{\boldsymbol{\gamma}}_{\!p}(\lambda)=(\widehat{\gamma}_0(\lambda), \dots, \widehat{\gamma}_p(\lambda))^\top\!,$ see  \eqref{Estimate gamma_h uni-variate}, and $\pmb{\mathbb{B}}(\lambda)\coloneqq (\mathbb{B}_0(\lambda), \dots, \mathbb{B}_p(\lambda))^{\!\top}\!,$ where $\mathbb{B}_0, \dots, \mathbb{B}_p$ are independent, standard Brownian motions on $[0,1].$ The proof is performed in several steps.

\paragraph{Step 1.} First, to verify the claim in Theorem \ref{piv_M_cor}, we prove that
\begin{align}\label{eq: process cvg gammas}
\sqrt{N}(\widehat{\boldsymbol{g}}_p - \boldsymbol{g}_p)\, \rightsquigarrow \,\pmb{\mathcal{G}} \coloneqq \big\{\Sigma_p^{1/2}\,\pmb{\mathbb{B}}(\lambda)\big\}_{\!\lambda \in [0,1]},
\end{align}
where $\Sigma_p$ is given in \eqref{AN for vector of autocov}. To establish this, observe that the innovations of the linear process are i.i.d.\ with finite fourth moments, and the coefficients satisfy $\sum_{j=1}^\infty j|\theta_j| < \infty$. By \citet[][Proposition~2.1]{HoermannKokoszka2010}, $(X_k)$ is $L^4$\nobreakdash-$m$\nobreakdash-approximable (see Remark \ref{rem:L^p-m}). It then follows by, e.g., \citet[][Lemma~B.1]{Kuehnert2022} that $(\boldsymbol{Y}_{\!\!k})_{k \in \mathbb{Z}}$, with $\boldsymbol{Y}_{\!\!k} \coloneqq (X_kX_{k+0}-\gamma_0, \dots, X_kX_{k+p}-\gamma_p)^\top \in \mathbb{R}^{p+1}$, is $L^2$\nobreakdash-$m$\nobreakdash-approximable, and by \citet[][Theorem~1.1]{Jirak2013},
\begin{align*}
    \bigg\{\frac{1}{\sqrt{N}}\!\sum_{k=1}^{\lfloor\lambda N\rfloor}\boldsymbol{Y}_{\!\!k}\bigg\}_{\!\lambda\in[0,1]}
    \rightsquigarrow \,\pmb{\mathcal{G}}.
\end{align*}

In what follows, let $m_N(\lambda) \coloneqq \max\{\lfloor\lambda N\rfloor-p,0\}$ for $\lambda\in[0,1]$. We first show that replacing $\lfloor\lambda N\rfloor$ by $m_N(\lambda)$ does not affect the limiting distribution. To this end, let $Z_k\coloneqq\max_{0\leq h\leq p}|X_kX_{k+h}|$. By stationarity, the Cauchy--Schwarz inequality, and the assumption of finite fourth moments, we have $\Erw(Z_0^2) \leq (p+1)\Erw(X_0^4) < \infty.$ Thus, for every $\epsilon>0$, the union bound yields
\begin{align*}
    \Prob\!\Big(\max_{1\leq k\leq N}Z_k>\epsilon\sqrt{N}\Big)
    \,\leq\, 
    N\Prob\!\Big(Z_0>\epsilon\sqrt{N}\Big)
    \,\leq\, \frac{1}{\epsilon^2}\Erw\Big[Z_0^2\mathbf{1}_{\{Z_0>\epsilon\sqrt{N}\}}\Big]
    \,\longrightarrow\,
    0,
\end{align*}
and thus $\max_{1\leq k\leq N}Z_k=o_{\Prob}(\sqrt{N})$. Since $0\leq\lfloor\lambda N\rfloor-m_N(\lambda)\leq p$, we obtain
\begin{align*}
\sup_{\lambda\in[0,1]}\max_{0\leq h\leq p}\bigg|\sum_{k=m_N(\lambda)+1}^{\lfloor\lambda N\rfloor}\boldsymbol{Y}_{\!\!k}\,\bigg|
\,\leq\,
p\Big(\max_{1\leq k\leq N}Z_k+\gamma_0\Big)
\,=\,
o_{\Prob}(\sqrt{N}).
\end{align*}
Therefore, by Slutsky's theorem, with $\pmb{\mathcal{G}}$ as in \eqref{eq: process cvg gammas},
\begin{align}\label{Prelim result process cvgc}
    \bigg\{\frac{1}{\sqrt{N}}\!\sum_{k=1}^{m_N(\lambda)}\boldsymbol{Y}_{\!\!k}\bigg\}_{\!\lambda\in[0,1]}
    \rightsquigarrow \,\pmb{\mathcal{G}}\,.
\end{align}
Further, $m_N(\lambda)\leq\lfloor\lambda(N-h)\rfloor$ for all $\lambda\in[0,1]$ and $h=0,\ldots,p$. Thus, the definition of $\widehat{\gamma}_h(\lambda)$ in \eqref{Estimate gamma_h uni-variate} gives for each component of the sum in \eqref{Prelim result process cvgc}:
\begin{align*}
    \sum_{k=1}^{m_N(\lambda)}\!Y_{k,h}
    \,=
    \sum_{k=1}^{m_N(\lambda)}\!X_kX_{k+h}-m_N(\lambda)\gamma_h
    \,=\,
    N\big(\widehat{\gamma}_h(\lambda)-\lambda\gamma_h\big)+S_h(\lambda),
\end{align*}
where
\[
    S_h(\lambda)\coloneqq\big(\lambda N-m_N(\lambda)\big)\gamma_h
    -\sum_{k=m_N(\lambda)+1}^{\lfloor\lambda(N-h)\rfloor}X_kX_{k+h}.
\]
Further, as $|\gamma_h|\leq\Erw\!|X_kX_{k+h}|\leq\gamma_0,$ and $0\leq\lambda N-m_N(\lambda)\leq p+1$ and $0\leq\lfloor\lambda(N-h)\rfloor-m_N(\lambda)\leq p,$ the sum defining $S_h(\lambda)$ contains at most $p$ terms. Consequently,
\[
|S_h(\lambda)|\leq(p+1)\gamma_0+p\max_{1\leq k\leq N}Z_k
\]
uniformly in $\lambda\in[0,1]$ and $h=0,\ldots,p$. Using the maximal bound established above, we obtain
\begin{align*}
    \big\|\big\{\big(S_0(\lambda),\dots,S_p(\lambda)\big)^\top\big\}_{\lambda\in[0,1]}\big\|_\infty
    \,=
    \sup_{\lambda\in[0,1]}\sup_{0\leq h\leq p}|S_h(\lambda)|
    \,=\,o_{\Prob}(\sqrt{N}).
\end{align*}

\noindent Together with \eqref{Prelim result process cvgc}, Slutsky's theorem, and the definitions of $\widehat{\boldsymbol{g}}_p,\boldsymbol{g}_p$ in \eqref{Def hatg_p and g_p}, we arrive at \eqref{eq: process cvg gammas}.

\paragraph{Step 2.} We now introduce the following map and establish its Hadamard differentiability:
\begin{align}\label{det87}
\phi\colon \!\!
\begin{cases}
    \,{\cal D}^+_{0,p+1} \to \ell^{\infty,p+1}([0,1]),\\
    \,\boldsymbol{f}_{\!p} \mapsto
    \phi(\boldsymbol{f}_{\!p})\colon \!\!
    \begin{cases}
        \,[0,1] \to \R^{p+1},\\
        \,\lambda \mapsto
        (\phi(\boldsymbol{f}_{\!p}))(\lambda)
        \!\coloneqq\!
        \begin{cases}
            \,\Big(
            \textstyle
            \frac{\det((f_{|j-i|}(\lambda))^0_{i,j=0})}
                 {\det((f_{|j-i|}(\lambda))^{-1}_{i,j=0})},
            \dots,
            \frac{\det((f_{|j-i|}(\lambda))^p_{i,j=0})}
                 {\det((f_{|j-i|}(\lambda))^{p-1}_{i,j=0})}
            \Big)^{\!\top}\!,
            & \!\lambda\in(0,1],\\
            \,\boldsymbol{0}, & \!\lambda=0,
        \end{cases}
    \end{cases}
\end{cases}
\end{align}
where $\det((f_{|j-i|}(\lambda))^{-1}_{i,j=0})=1$ by convention and where its domain is defined as
\begin{align*}
{\cal D}^+_{0,p+1} \coloneqq
\Big\{
\boldsymbol{f}_{\!p}\in
\ell^{\infty,p+1}_0([0,1])
\,\Big|\,
(f_{|j-i|}(\lambda))^p_{i,j=0}\succ0
\textnormal{ for all }\lambda\in(0,1]
\Big\},
\end{align*}
with $\succ0$ denoting positive definiteness. Here, $\boldsymbol{f}_{\!p} =\{\boldsymbol{f}_{\!p}(\lambda)\}_{\lambda\in[0,1]} =\{(f_0(\lambda),\dots,f_p(\lambda))^{\!\top}\}_{\lambda\in[0,1]}$. Since every leading principal submatrix of a positive definite matrix is positive definite, all determinants in \eqref{det87} are positive for
$\lambda\in(0,1]$. Moreover, by the Schur complement formula in \eqref{Schur}
and the positive definiteness of the corresponding
Toeplitz matrices, we have
\begin{align*}
0 <
\frac{
\det((f_{|j-i|}(\lambda))^k_{i,j=0})
}{
\det((f_{|j-i|}(\lambda))^{k-1}_{i,j=0})
}
\le f_0(\lambda)
\le \|f_0\|_\infty < \infty
\end{align*}
for all $\lambda\in(0,1]$ and $k=0,\dots,p.$ Hence, $\phi$ is well-defined as a map into $\ell^{\infty,p+1}([0,1]).$ Moreover, by its definition, with $\boldsymbol{g}_p$ in \eqref{Def hatg_p and g_p} and $M_p$ in \eqref{Def M_p}, we obtain
\begin{align}\label{phi_p at g_p}
    (\phi(\boldsymbol{g}_p))(\lambda)
    =\lambda(M_0,\dots,M_p)^\top,
    \qquad \lambda\in[0,1].
\end{align}

To establish the Hadamard differentiability of $\phi$,
we first consider its component maps $\phi^{(k)}$,
$k=1,\dots,p+1$, defined on the common domain
${\cal D}^+_{0,p+1}$ by
\begin{align*}
\phi^{(k)}\colon \!
\begin{cases}
    \;{\cal D}^+_{0,p+1}\rightarrow\ell^\infty([0,1]),\\
    \;\boldsymbol{f}_{\!p}\mapsto
    \phi^{(k)}(\boldsymbol{f}_{\!p})\colon\!
    \begin{cases}
        \;[0,1]\rightarrow\R,\\
        \;\lambda\mapsto
        \big(\phi^{(k)}(\boldsymbol{f}_{\!p})\big)(\lambda)
        \coloneqq
        \begin{cases}
            \;\textstyle
            \frac{\det((f_{|j-i|}(\lambda))^{k-1}_{i,j=0})}
                 {\det((f_{|j-i|}(\lambda))^{k-2}_{i,j=0})},
            & \!\lambda\in(0,1],\\
            \;0, & \!\lambda=0.
        \end{cases}
    \end{cases}
\end{cases}
\end{align*}

\noindent By the definitions of $\phi^{(k)}$ and $\boldsymbol{g}_p$
in \eqref{Def hatg_p and g_p}, we have
\begin{align*}
    (\phi^{(k)}(\boldsymbol{g}_p))(\lambda)
    =\lambda M_{k-1},
    \qquad \lambda\in[0,1],\quad k=1,\dots,p+1.
\end{align*}
Furthermore, for each $k=1,\dots,p+1$, consider the map
\begin{align}\label{Def_tilde_phik}
    \widetilde{\phi}^{(k)}\colon
    \;\mathcal U_{p+1}\to\mathbb R, ~ \boldsymbol f_{\!p} = (f_0, \dots, f_p)^\top \mapsto
    \dfrac{\det((f_{|j-i|})^{k-1}_{i,j=0})}
    {\det((f_{|j-i|})^{k-2}_{i,j=0})},
\end{align}
where
\begin{align*}
    \mathcal U_{p+1}\coloneqq
    \big\{
    \boldsymbol f_{\!p} = (f_0, \dots, f_p)^\top\in\mathbb R^{p+1}
    \,\big|\,
    (f_{|j-i|})^p_{i,j=0}\succ0
    \big\}.
\end{align*}
Since $\mathcal U_{p+1}$ is open, the determinants are polynomial functions, and the denominators are strictly positive on $\mathcal U_{p+1}$, it follows from the quotient rule that $\widetilde{\phi}^{(k)}$ is continuously differentiable on $\mathcal U_{p+1}$. By assumption,
$\boldsymbol g_p(\lambda)
=\lambda\boldsymbol\gamma_p\in\mathcal U_{p+1}$
for every $\lambda\in(0,1]$.
Further, by using \eqref{grad_identity}, the gradient of $\widetilde{\phi}^{(k)}$ evaluated at $\boldsymbol g_p(\lambda)$ satisfies
\begin{align*}
\nabla\widetilde{\phi}^{(k)}
(\boldsymbol g_p(\lambda))
=
\nabla M_{k-1,\boldsymbol\gamma_p},
\qquad \lambda\in(0,1].
\end{align*}
Consequently, for each fixed $\lambda\in(0,1]$,
\begin{align*}
    &\Big|
    \widetilde{\phi}^{(k)}
    (\boldsymbol g_p(\lambda)+\boldsymbol z_p)
    -\widetilde{\phi}^{(k)}
    (\boldsymbol g_p(\lambda)) - \nabla M_{k-1,\boldsymbol\gamma_p}^{\top}
    \boldsymbol z_p
    \Big|
    =o(\|\boldsymbol z_p\|),
    \qquad \|\boldsymbol z_p\|\to0,
\end{align*}
provided that $\boldsymbol g_p(\lambda)+\boldsymbol z_p \in \mathcal{U},$ where $\|\cdot\|$ denotes a norm on
$\mathbb R^{p+1}.$ This motivates the definition of the bounded
linear operators for $k=1, \dots, p+1,$
\begin{align}\label{det1}
    \phi^{(k)\prime}_{\boldsymbol{g}_p}\colon
    \begin{cases}
    \;\mathcal C_{0,p+1}
    \rightarrow\ell^\infty([0,1]),\\
    \;\boldsymbol{f}_{\!p}\mapsto
    \phi^{(k)\prime}_{\boldsymbol{g}_p}
    (\boldsymbol{f}_{\!p})\colon
    \begin{cases}
    \;[0,1]\to\mathbb R,\\
    \;\lambda\mapsto
    \big(\phi^{(k)\prime}_{\boldsymbol{g}_p}
    (\boldsymbol{f}_{\!p})\big)(\lambda)
    \coloneqq
    \nabla M^\top_{k-1,\boldsymbol{\gamma}_p}
    \boldsymbol{f}_{\!p}(\lambda),
    \end{cases}
    \end{cases}
\end{align}
where, in this definition,
\begin{align}\label{setC0p+1}
    \mathcal C_{0,p+1}
    \coloneqq
    \big\{
    \boldsymbol h_p\in C([0,1],\mathbb R^{p+1})
    \,\big|\,
    \boldsymbol h_p(0)=\boldsymbol{0}
    \big\},
\end{align}
equipped with the supremum norm, and $C([0,1],\mathbb R^{p+1})$ denotes the space of continuous functions from $[0,1]$ to $\mathbb{R}^{p+1}.$ We proceed by establishing Hadamard differentiability of every $\phi^{(k)}$ at $\boldsymbol g_p$, tangentially to the space $\mathcal C_{0,p+1} \subset \ell^{\infty,p+1}_0$ \citep[see][pp.\ 296--297, for a formal definition]{vanderVaart1998}. Let $\boldsymbol h_p\in\mathcal C_{0,p+1}$ and let $\boldsymbol h_p^t\in
\ell^{\infty,p+1}_0([0,1])$
satisfy
$\|\boldsymbol h_p^t-\boldsymbol h_p\|_\infty
\to0$ as $t\downarrow0$, with
$\boldsymbol g_p+t\boldsymbol h_p^t
\in{\cal D}^+_{0,p+1}$
for sufficiently small $t>0$. Since $\boldsymbol g_p(0)
=\boldsymbol h_p(0)
=\boldsymbol h_p^t(0)=\boldsymbol{0}$, we have
\begin{align*}
&\bigg\|
\frac{
\phi^{(k)}(\boldsymbol g_p+t\boldsymbol h_p^t)
-\phi^{(k)}(\boldsymbol g_p)
}{t}
-\phi^{(k)\prime}_{\boldsymbol g_p}
(\boldsymbol h_p)
\bigg\|_\infty\\
&\quad=
\sup_{\lambda\in(0,1]}
\bigg|
\frac{
\widetilde{\phi}^{(k)}
(\boldsymbol g_p(\lambda)
+t\boldsymbol h_p^t(\lambda))
-\widetilde{\phi}^{(k)}
(\boldsymbol g_p(\lambda))
}{t} -\nabla M_{k-1,\boldsymbol\gamma_p}^{\top}
\boldsymbol h_p(\lambda)
\bigg|.
\end{align*}

In the following, we fix $\epsilon\in(0,1)$ and consider the suprema over the subintervals $(0,\epsilon]$ and $[\epsilon,1]$ separately, starting with the latter. Since the map $\lambda\mapsto\boldsymbol g_p(\lambda)=\lambda\boldsymbol\gamma_p$ is continuous, and the Toeplitz matrix associated with $\boldsymbol\gamma_p$ is positive definite by assumption, $K_\epsilon \coloneqq
\{\boldsymbol g_p(\lambda):\lambda\in[\epsilon,1]\} \subset\mathcal U_{p+1}$. Further, as $\mathcal U_{p+1}$ is open and $K_\epsilon$ is compact, there exists $\delta>0$ such that the closed $\delta$-neighborhood of $K_\epsilon$ is compact and contained in $\mathcal U_{p+1}$. As the map $\widetilde{\phi}^{(k)}$ defined in \eqref{Def_tilde_phik} is continuously differentiable on $\mathcal U_{p+1}$, its gradient is continuous and therefore uniformly continuous on this compact neighborhood. Consequently, the first-order Taylor expansion holds uniformly for $\lambda\in[\epsilon,1]$. Together with the uniform convergence $\boldsymbol h_p^t\to\boldsymbol h_p$, this yields
\begin{align*}
    &\sup_{\lambda\in[\epsilon,1]}
    \bigg|
    \frac{
    \widetilde{\phi}^{(k)}
    (\boldsymbol g_p(\lambda)
    +t\boldsymbol h_p^t(\lambda))
    -\widetilde{\phi}^{(k)}
    (\boldsymbol g_p(\lambda))
    }{t} - \nabla M_{k-1,\boldsymbol\gamma_p}^{\top}
    \boldsymbol h_p(\lambda)
    \bigg|
    \stackrel{t\downarrow0}{\longrightarrow} 0.
\end{align*}
Next, consider the subset $(0,\epsilon] \subset (0,1].$ We show that each map $\widetilde{\phi}^{(k)}$  with $k=1, \dots, p+1,$ is Lipschitz continuous on the entire set $\mathcal U_{p+1},$ with Lipschitz constant $L_k = 4^{k-1}.$ Hereto, let $\boldsymbol f_p = (f_0, \dots, f_p)^\top\!,$ $\boldsymbol f'_p = (f'_0, \dots, f'_p)^\top \in \mathcal U_{p+1},$ and define
\[
    \boldsymbol a_{k-1}
    \coloneqq
    T^{-1}_{k-2}(\boldsymbol f_p)
    \boldsymbol r_{k-1}(\boldsymbol f_p) \in \mathbb{R}^{k-1},
    \qquad
    \boldsymbol a'_{k-1}
    \coloneqq
    T^{-1}_{k-2}(\boldsymbol f_p')
    \boldsymbol r_{k-1}(\boldsymbol f_p') \in \mathbb{R}^{k-1}, \qquad k=2, \dots, p+1,
\]
where $T_{k-2}(\boldsymbol f_p) \coloneqq (f_{|j-i|})_{i,j=0}^{k-2} \in \mathbb{R}^{(k-1)\times (k-1)}$ and $\boldsymbol r_{k-1}(\boldsymbol f_p) \coloneqq (f_{k-1},\ldots,f_1)^\top \in \mathbb{R}^{k-1}.$ By the definition of $\widetilde{\phi}^{(k)}$ in \eqref{Def_tilde_phik}, we have
\[
\widetilde{\phi}^{(1)}(\boldsymbol f_p)
-\widetilde{\phi}^{(1)}(\boldsymbol f_p')
=f_0-f_0',
\]
which immediately establishes Lipschitz continuity of $\widetilde{\phi}^{(1)}$ with constant $L_1=1$. Further, for $k=2,\ldots,p+1$, the Schur complement formula in \eqref{Schur} yields
\begin{align*}
    \widetilde{\phi}^{(k)}(\boldsymbol f_p)
    -\widetilde{\phi}^{(k)}(\boldsymbol f_p')
    &=
    f_0-f_0'
    -\boldsymbol r_{k-1}(\boldsymbol f_p)^\top
    \boldsymbol a_{k-1}
    +\boldsymbol r_{k-1}(\boldsymbol f_p')^\top
    \boldsymbol a_{k-1}'\\
    &=
    \boldsymbol v_k^\top
    T_{k-1}(\boldsymbol f_p)\boldsymbol v_k
    -(\boldsymbol v_k')^\top
    T_{k-1}(\boldsymbol f_p')\boldsymbol v_k',
\end{align*}
where
\[
\boldsymbol v_k
\coloneqq
\big(-\boldsymbol a_{k-1}^{\top},1\big)^{\top}
\in\mathbb R^k,
\qquad
\boldsymbol v_k'
\coloneqq
\big(-(\boldsymbol a_{k-1}')^{\top},1\big)^{\top}
\in\mathbb R^k.
\]
Since $T_{k-2}(\boldsymbol f_p)$ and
$T_{k-2}(\boldsymbol f_p')$ are positive definite,
the vectors $\boldsymbol v_k$ and $\boldsymbol v_k'$
minimize the respective quadratic forms among all
vectors in $\mathbb R^k$ whose last component equals
one. In particular,
\begin{align*}
    \boldsymbol v_k^\top
    \big(T_{k-1}(\boldsymbol f_p)-T_{k-1}(\boldsymbol f_p')\big)
    \boldsymbol v_k
    \,\leq\,
    \widetilde{\phi}^{(k)}(\boldsymbol f_p)
    -\widetilde{\phi}^{(k)}(\boldsymbol f_p')
    \,\leq\,
    (\boldsymbol v_k')^\top
    \big(T_{k-1}(\boldsymbol f_p)-T_{k-1}(\boldsymbol f_p')\big)
    \boldsymbol v_k'.
\end{align*}
Furthermore, positive definiteness and the Levinson--Durbin recursion imply that the $1$-norms of the corresponding vectors increase by at most a factor of two at each step. Hence, for the $1$-norms holds $\|\boldsymbol v_k\|_1 \leq 2^{k-1}$ and $\|\boldsymbol v_k'\|_1\leq 2^{k-1}.$ Combining these bounds, using the Toeplitz structure and the triangle inequality,
we obtain
\begin{align*}
    \left|
    \widetilde{\phi}^{(k)}(\boldsymbol f_p)
    -\widetilde{\phi}^{(k)}(\boldsymbol f_p')
    \right|
    &\leq
    \max_{\boldsymbol w_k\in\{\boldsymbol v_k,\boldsymbol v_k'\}}
    \left|
    \boldsymbol w_k^\top
    \big(
    T_{k-1}(\boldsymbol f_p)-T_{k-1}(\boldsymbol f_p')
    \big)
    \boldsymbol w_k
    \right|\\
    &\leq
    \max_{\boldsymbol w_k\in\{\boldsymbol v_k,\boldsymbol v_k'\}}
    \|\boldsymbol w_k\|_1^2
    \|\boldsymbol f_p-\boldsymbol f_p'\|_\infty\\
    &\leq
    4^{k-1}
    \|\boldsymbol f_p-\boldsymbol f_p'\|_\infty.
\end{align*}
\noindent 
Thus, as claimed, every $\widetilde{\phi}^{(k)},$ with $k=1, \dots, p+1,$ is Lipschitz
continuous on $\mathcal U_{p+1}$ with Lipschitz
constant $L_k=4^{k-1}.$ Therefore,
\begin{align*}
    \sup_{\lambda\in(0,\epsilon]}
    \bigg|
    \frac{
    \widetilde{\phi}^{(k)}
    (\boldsymbol g_p(\lambda)
    +t\boldsymbol h_p^t(\lambda))
    -\widetilde{\phi}^{(k)}
    (\boldsymbol g_p(\lambda))
    }{t} 
    -\nabla M_{k-1,\boldsymbol\gamma_p}^{\top}
    \boldsymbol h_p(\lambda)
    \bigg|
    &\,\le
    L_k\|\boldsymbol h_p^t-\boldsymbol h_p\|_\infty\\
    &\,\quad\;+
    \big(
    L_k+\|\nabla M_{k-1,\boldsymbol\gamma_p}\|
    \big)\!
    \sup_{\lambda\in[0,\epsilon]}
    \|\boldsymbol h_p(\lambda)\|.
\end{align*}
The first term
converges to zero, and the second term can be made arbitrarily small by choosing $\epsilon>0$ sufficiently small, since $\boldsymbol h_p$ is continuous at zero with $\boldsymbol h_p(0)=0.$ Therefore, altogether,
\begin{align*}
    &\bigg\|
    \frac{
    \phi^{(k)}(\boldsymbol g_p+t\boldsymbol h_p^t)
    -\phi^{(k)}(\boldsymbol g_p)
    }{t}
    -\phi^{(k)\prime}_{\boldsymbol g_p}
    (\boldsymbol h_p)
    \bigg\|_\infty \stackrel{t\downarrow0}{\longrightarrow} 0,
\end{align*}
so $\phi^{(k)} : {\cal D}^+_{0,p+1} \to \ell^{\infty}([0,1])$ is Hadamard differentiable at $\boldsymbol g_p$, tangentially to the space $\mathcal C_{0,p+1} \subset \ell^{\infty,p+1}_0([0,1])$, with derivative $\phi^{(k)\prime}_{\boldsymbol g_p}$ in \eqref{det1}.

Consequently, the vector-valued map
$\phi\colon {\cal D}^+_{0,p+1}\subset
\ell^{\infty,p+1}_0([0,1])
\to\ell^{\infty,p+1}([0,1])$
is Hadamard differentiable at $\boldsymbol g_p\in{\cal D}^+_{0,p+1}$,
tangential to the space $\mathcal C_{0,p+1} \subset \ell^{\infty,p+1}_0([0,1])$,
with derivative
\begin{align*}
    \phi'_{\!\boldsymbol{g}_p}\colon
    \begin{cases}
    \;\mathcal C_{0,p+1}
    \rightarrow\ell^{\infty,p+1}([0,1]),\\
    \;\boldsymbol{f}_{\!p}\mapsto
    \phi'_{\!\boldsymbol{g}_p}(\boldsymbol{f}_{\!p})
    \colon
    \begin{cases}
    \;[0,1]\to\R^{p+1},\\
    \;\lambda\mapsto
    \big(\phi'_{\!\boldsymbol{g}_p}
    (\boldsymbol{f}_{\!p})\big)(\lambda)
    \coloneqq
    \mathcal{M}_{p,\boldsymbol{\gamma}_p}
    \boldsymbol{f}_{\!p}(\lambda),
    \end{cases}
    \end{cases}
\end{align*}
where
\begin{align}
    \mathcal{M}_{p,\boldsymbol{\gamma}_p}
    \coloneqq
    \begin{pmatrix}
    \nabla M^\top_{0,\boldsymbol{\gamma}_p}\\
    \nabla M^\top_{1,\boldsymbol{\gamma}_p}\\
    \vdots\\
    \nabla M^\top_{p,\boldsymbol{\gamma}_p}
    \end{pmatrix}
    =
    \begin{pmatrix}
    \partial_{\gamma_0}M_0(\boldsymbol{\gamma}_p)
    & 0 & \cdots & 0\\
    \partial_{\gamma_0}M_1(\boldsymbol{\gamma}_p)
    & \partial_{\gamma_1}M_1(\boldsymbol{\gamma}_p)
    & \ddots & \vdots\\
    \vdots & \vdots & \ddots & 0\\
    \partial_{\gamma_0}M_p(\boldsymbol{\gamma}_p)
    & \partial_{\gamma_1}M_p(\boldsymbol{\gamma}_p)
    & \cdots &
    \partial_{\gamma_p}M_p(\boldsymbol{\gamma}_p)
    \end{pmatrix}.
    \label{kue7}
\end{align}

\paragraph{Step 3.}
Since the limiting process $\pmb{\mathcal G}$ in
\eqref{eq: process cvg gammas} has continuous sample
paths vanishing at zero, it belongs to
$\mathcal C_{0,p+1}$ defined in \eqref{setC0p+1}
almost surely. Thus, the Hadamard differentiability
established in Step 2, tangentially to
$\mathcal C_{0,p+1}$, would allow us to apply the
functional delta method
\citep[][Theorem~20.8]{vanderVaart1998}
if $\widehat{\boldsymbol g}_p$ belonged to
$\mathcal D^+_{0,p+1}$, the domain of $\phi$ in
\eqref{det87}. However, this condition is not
satisfied, since
$\widehat{\boldsymbol\gamma}_p(\lambda)
=\boldsymbol 0$ for all $\lambda\in(0,1/N)$,
so that the corresponding Toeplitz matrices are
singular. To overcome this difficulty, we construct
an auxiliary process taking values in
$\mathcal D^+_{0,p+1}$ and show that it is
asymptotically equivalent to the original
empirical process.

To this end, define
\begin{align*}
    c_N(\lambda)
    &\coloneqq
    \big[-\lambda_{\min}(\widehat G_p(\lambda))\big]_+
    +\frac{\lambda}{N},
\end{align*}
where $[x]_+=\max\{x,0\},$ $\lambda_{\min}$ refers to the smallest eigenvalue, and $\widehat G_p(\lambda) = (\widehat\gamma_{j-i}(\lambda))_{i,j=0}^p.$ Further,
\begin{align}\label{auxiliary_g_process}
    \widetilde{\boldsymbol g}_p(\lambda)
    \coloneqq
    \widehat{\boldsymbol\gamma}_p(\lambda)
    +c_N(\lambda)\boldsymbol e_1,
    \qquad \lambda\in[0,1],
\end{align}
where
$\boldsymbol e_1=(1,0,\ldots,0)^\top\in\mathbb R^{p+1}$.
The associated Toeplitz matrices satisfy
\begin{align*}
    \widetilde G_{p,N}(\lambda)
    =
    \widehat G_p(\lambda)+c_N(\lambda)I_{p+1},
    \qquad
    \lambda_{\min}(\widetilde G_{p,N}(\lambda))
    \geq \frac{\lambda}{N}>0,
    \qquad 
    \lambda\in(0,1].
\end{align*}
Moreover, $c_N(0)=0$. Hence, for the auxiliary process holds $\widetilde{\boldsymbol g}_p = \{\widetilde{\boldsymbol g}_p(\lambda)\}_{\lambda \in [0,1]}
\in\mathcal D^+_{0,p+1}$.

Next, we establish the asymptotic equivalence
of the two processes. In what follows, let $\epsilon \in (0,1).$ Since the population Toeplitz matrix $G_p$ is positive definite, we have
\[
    \inf_{\lambda\in[\epsilon,1]}
    \lambda_{\min}(\lambda G_p)
    =
    \epsilon\lambda_{\min}(G_p)>0.
\]
Hence, by Weyl's inequality,
\[
    \inf_{\lambda\in[\epsilon,1]}\lambda_{\min}(\widehat G_p(\lambda))
    \geq
    \epsilon\lambda_{\min}(G_p) - \sup_{\lambda\in[0,1]}\|\widehat G_p(\lambda)-\lambda G_p\|_{\mathrm{op}}\,,
\]
with $\|\cdot\|_{\mathrm{op}}$ denoting the operator norm. Therefore, since \eqref{eq: process cvg gammas} implies $\sup_{\lambda\in[0,1]}
\|\widehat G_p(\lambda)-\lambda G_p\|_{\mathrm{op}}
=O_{\Prob}(N^{-1/2}),$ we obtain
\[
    \Prob\!\Big(
    \widehat G_p(\lambda)\succ0
    \text{ for all }\lambda\in[\epsilon,1]
    \Big) \longrightarrow 1.
\]

For $\lambda\in[0,\epsilon]$, Weyl's inequality and
$\lambda_{\min}(\lambda G_p)\geq0$ imply
\begin{align*}
    \big[-\lambda_{\min}(
    \widehat G_p(\lambda))\big]_+
    &\leq
    \|\widehat G_p(\lambda)
    -\lambda G_p\|_{\mathrm{op}}.
\end{align*}
Hence, by the Toeplitz structure and with $\boldsymbol Z_N \coloneqq \sqrt N(\widehat{\boldsymbol g}_p - \boldsymbol g_p),$ one obtains
\begin{align*}
    \sqrt N\,c_N(\lambda)
    &\leq
    \sqrt N
    \|\widehat G_p(\lambda)
    -\lambda G_p\|_{\mathrm{op}}
    +\frac{\lambda}{\sqrt N}\\
    &\leq
    (p+1)\|\boldsymbol Z_N(\lambda)\|_\infty
    +\frac{1}{\sqrt N}.
\end{align*}
Further, since $\boldsymbol Z_N\rightsquigarrow\pmb{\mathcal G}$ by \eqref{eq: process cvg gammas}, the continuous mapping theorem and the Portmanteau theorem applied to the closed set $[\eta,\infty)$ imply that, for every fixed $\epsilon>0$ and $\eta>0$,
\[
    \limsup_{N\to\infty}
    \Prob\Big(
    \sup_{\lambda\in[0,\epsilon]}
    \|\boldsymbol Z_N(\lambda)\|_\infty>\eta
    \Big) 
    \leq
    \Prob\Big(
    \sup_{\lambda\in[0,\epsilon]}
    \|\pmb{\mathcal G}(\lambda)\|_\infty\geq\eta
    \Big).
\]
Moreover, $\pmb{\mathcal G}$ has continuous sample paths vanishing at zero,
\[
    \lim_{\epsilon\downarrow0}
    \limsup_{N\to\infty}
    \Prob\Big(
    \sup_{\lambda\in[0,\epsilon]}
    \|\boldsymbol Z_N(\lambda)\|_\infty>\eta
    \Big)=0.
\]

Consequently, the contribution from
$[0,\epsilon]$ can be made arbitrarily small
in probability by choosing $\epsilon>0$
sufficiently small. On the other hand,
the positive definiteness established
uniformly on $[\epsilon,1]$ implies that
$c_N(\lambda)=\lambda/N$ on this interval
with probability tending to one.
Combining both bounds yields
\begin{align}\label{auxiliary_equivalence}
    \sqrt N
    \|\widetilde{\boldsymbol g}_p
    -\widehat{\boldsymbol g}_p\|_\infty
    &=
    \sqrt N\sup_{\lambda\in[0,1]}c_N(\lambda)
    \xrightarrow{\Prob}0.
\end{align}

\paragraph{Step 4.} First, \eqref{eq: process cvg gammas}, \eqref{auxiliary_equivalence} and Slutsky's theorem yield for the auxiliary
process,
\begin{align*}
    \sqrt N
    (\widetilde{\boldsymbol g}_p
    -\boldsymbol g_p)
    \rightsquigarrow\pmb{\mathcal G}.
\end{align*}
As $\widetilde{\boldsymbol g}_p
\in\mathcal D^+_{0,p+1}$, we are allowed to apply the functional delta method. Together with \eqref{phi_p at g_p}, this yields
\begin{align}\label{auxiliarydeltalimit}
    \sqrt N
    \big(
    \phi(\widetilde{\boldsymbol g}_p)
    -\phi(\boldsymbol g_p)
    \big)
    \rightsquigarrow
    \phi'_{\boldsymbol g_p}(\pmb{\mathcal G})
    =
    \big\{
    \mathcal M_{p,\boldsymbol\gamma_p}
    \Sigma_p^{1/2}\,
    \pmb{\mathbb B}(\lambda)
    \big\}_{\lambda\in[0,1]},
\end{align}
where the limit agrees with the one in our claim. 

It remains to establish the asymptotic equivalence
of the transformed auxiliary process and the
original sequential estimators. To this end, let
\begin{align*}
    A_N
    \coloneqq
    \big\{\lambda\in(0,1]:
    \widehat G_p(\lambda)\succ0\big\},
    \qquad 
    B_N
    \coloneqq
    \big\{\lambda\in(0,1]:
    \widehat G_p(\lambda)\not\succ0\big\}.
\end{align*}
On $A_N$, Lipschitz continuity of the component-wise functions established in Step 2 implies, with $\widehat{\boldsymbol M}_p \coloneqq (\widehat{M}_0(\lambda),  \dots, \widehat{M}_p(\lambda))^{\!\top}\!,$ and with $\widehat{M}_k(\lambda)$ for each $k$ defined as in \eqref{Def Mp gamma},
\begin{align*}
    \sup_{\lambda\in A_N}
    \big\|
    \widehat{\boldsymbol M}_p(\lambda)
    -\phi(\widetilde{\boldsymbol g}_p)(\lambda)
    \big\|_\infty
    \leq
    L\sup_{\lambda\in[0,1]}c_N(\lambda),
\end{align*}
where
$L\coloneqq\max_{1\leq k\leq p+1}L_k<\infty$.
Hence, by \eqref{auxiliary_equivalence},
the right-hand side is $o_{\Prob}(N^{-1/2})$. Further, on the event $B_N$, positive definiteness of $G_p$
and Weyl's inequality give
\begin{align}\label{rellambdaZ}
    \sqrt N\,\lambda\lambda_{\min}(G_p)
    \,\leq\,
    \sqrt N
    \|\widehat G_p(\lambda)-\lambda G_p\|_{\mathrm{op}}
    \,\leq\,
    (p+1)\|\boldsymbol Z_N(\lambda)\|_\infty.
\end{align}
By the arguments in Step 3 and the bound \eqref{rellambdaZ}, it holds
\[
    \sup_{\lambda\in B_N}\|\boldsymbol Z_N(\lambda)\| = o_{\Prob}(1),
    \qquad
    \sup_{\lambda\in B_N}\lambda
    =o_{\Prob}(N^{-1/2}).
\]
By the Schur complement representation in Step 2
and the definition of the sequential estimators,
the components of
$\widehat{\boldsymbol M}_p(\lambda)$
and $\phi(\widetilde{\boldsymbol g}_p)(\lambda)$
are non-negative and bounded above by
$\widehat\gamma_0(\lambda)$ and
$\widetilde g_{0,N}(\lambda)$, respectively.
Hence, by using 
$\widetilde g_{0,N}(\lambda)
=\widehat\gamma_0(\lambda)+c_N(\lambda)$
and
$\widehat\gamma_0(\lambda)
=\lambda\gamma_0+N^{-1/2}Z_{N,0}(\lambda)$,
we obtain, together with
\eqref{auxiliary_equivalence},
\begin{align*}
    \sup_{\lambda\in B_N}\big\|
    \widehat{\boldsymbol M}_p(\lambda)
    -\phi(\widetilde{\boldsymbol g}_p)(\lambda)
    \big\|_\infty
    &\leq
    2\sup_{\lambda\in B_N}\widehat\gamma_0(\lambda) + \sup_{\lambda\in B_N}c_N(\lambda)\\
    &\leq
    2\gamma_0\sup_{\lambda\in B_N}\lambda
    +2N^{-1/2}\!\sup_{\lambda\in B_N}\!\|\boldsymbol Z_N(\lambda)\|_\infty
    +\sup_{\lambda\in B_N}c_N(\lambda)\\
    &=o_{\Prob}(N^{-1/2}).
\end{align*}
Together with the bound on $A_N$, this yields
\[
    \big\|
    \widehat{\boldsymbol M}_p
    -\phi(\widetilde{\boldsymbol g}_p)
    \big\|_\infty
    =o_{\Prob}(N^{-1/2}).
\]
Finally, combining the latter with \eqref{auxiliarydeltalimit} and Slutsky's theorem, proves the claim.\hfill\qed

\subsection{Proof of Corollary \ref{piv_M_cor}}
An application of the continuous mapping theorem to the process in \eqref{det94} yields
\begin{align*}
    \sqrt{N}\big\{
    \widehat{M}_p(\lambda)-\lambda M_p
    \big\}_{\!\lambda\in[0,1]}
    \rightsquigarrow
    \big\{
    \tau_p\mathbb{B}(\lambda)
    \big\}_{\!\lambda\in[0,1]}\,,
\end{align*}
where $\tau_p^2=(\nabla M_{p,\boldsymbol{\gamma}_p})^\top\Sigma_p\nabla M_{p,\boldsymbol{\gamma}_p}$. Since $\Sigma_p$ is positive definite and $\nabla M_{p,\boldsymbol{\gamma}_p}\neq0$ by Lemma \ref{LemA1}, we have $\tau_p>0$. Thus, another application of the continuous mapping theorem with the map
\[
    \ell^\infty([0,1])\ni
    f=\{f(\lambda)\}_{\lambda\in[0,1]}
    \mapsto
    \frac{f(1)}
    {\int_0^1
    |f(\lambda)-\lambda f(1)|\,\mathrm{d}\lambda}\,,
\]
which is continuous whenever the denominator is non-zero, proves the claim. \hfill\qed

\subsection{Proof of Theorem \ref{thm: Distr. free limit of M'}}\label{SecA3}

We define the map
\begin{align}\label{det101}
\phi\colon \!
\begin{cases}
\;{\cal D}^S_{0,p+1} \rightarrow \ell^{\infty}([0,1]),\\
\;\boldsymbol{f}_{\!p} \mapsto
\phi(\boldsymbol{f}_{\!p})\colon \!
\begin{cases}
\;[0,1]\rightarrow\R,\\
\;\lambda\mapsto
(\phi(\boldsymbol{f}_{\!p}))(\lambda)
\coloneqq
\begin{cases}
\;\lambda f_p(\lambda)/f_0(\lambda),
& f_0(\lambda)>0,\;\lambda\in(0,1],\\
\;0, & \text{otherwise},
\end{cases}
\end{cases}
\end{cases}
\end{align}
where
$\boldsymbol{f}_{\!p}
=\{\boldsymbol{f}_{\!p}(\lambda)\}_{\lambda\in[0,1]}
=\{(f_0(\lambda),\dots,f_p(\lambda))^\top\}_{\lambda\in[0,1]}$
and the domain is defined as
\begin{align}\label{kue1}
{\cal D}^S_{0,p+1}\coloneqq
\Big\{
\boldsymbol{f}_{\!p}\in
\ell^{\infty,p+1}_0([0,1])
\,\Big|\,
0\leq f_p(\lambda)\leq f_0(\lambda)
\text{ for all }\lambda\in[0,1]
\Big\}.
\end{align}
The expression ``$\leq 1$'' in the notation indicates that the quotient bounded above by 1.
In particular, $0\leq(\phi(\boldsymbol{f}_{\!p}))(\lambda)
\leq\lambda$ for all $\lambda\in[0,1]$, so
$\phi$ is well-defined as a map into
$\ell^\infty([0,1])$.

With $\boldsymbol{M}_{\!p}
\coloneqq(M_0,\dots,M_p)^\top\in\R^{p+1}$
and
$\boldsymbol m_p
\coloneqq\{\lambda\boldsymbol{M}_{\!p}\}_{\lambda\in[0,1]}
\in{\cal D}$,
we obtain
\begin{align}\label{phi_p at g_p: Cor Mp/M0}
(\phi(\boldsymbol m_p))(\lambda)
=\lambda S_p
=\lambda\frac{M_p}{M_0},
\qquad \lambda\in[0,1],
\end{align}
where $M_0=\gamma_0>0$.

Furthermore, consider the map
\begin{align*}
\widetilde{\phi}\colon
\mathcal V^S_{p+1}\rightarrow\R,
\qquad
\boldsymbol f_{\!p}=(f_0,\dots,f_p)^\top
\mapsto\frac{f_p}{f_0},
\end{align*}
where $\mathcal V^S_{p+1}
\coloneqq
\{
\boldsymbol f_{\!p}\in\R^{p+1}
\,|\,f_0>0\}.$
Since $\mathcal V^S_{p+1}$ is open,
$\widetilde{\phi}$ is continuously differentiable
on $\mathcal V^S_{p+1}$. Its gradient at
$\lambda\boldsymbol{M}_{\!p}$ satisfies
\begin{align}\label{eq:identity gradient fp/f0}
\nabla\widetilde{\phi}
(\lambda\boldsymbol{M}_{\!p})
=\frac{1}{\lambda\gamma_0}s_p^\top,
\qquad \lambda\in(0,1],
\end{align}
where
$s_p\coloneqq(-S_p,0,\dots,0,1)$. Next, consider the bounded linear operator
\begin{align*}
\phi'_{\boldsymbol m_p}\colon
\begin{cases}
\;\mathcal C_{0,p+1}
\rightarrow\ell^\infty([0,1]),\\
\;\boldsymbol h_p\mapsto
\phi'_{\boldsymbol m_p}(\boldsymbol h_p)
\colon
\begin{cases}
\;[0,1]\rightarrow\R,\\
\;\lambda\mapsto
\big(\phi'_{\boldsymbol m_p}
(\boldsymbol h_p)\big)(\lambda)
\coloneqq
\frac{1}{\gamma_0}
s_p\boldsymbol h_p(\lambda),
\end{cases}
\end{cases}
\end{align*}
where
$\mathcal C_{0,p+1}
$ is defined as in \eqref{setC0p+1} and equipped
with the supremum norm. Analogous as in 
Section \ref{SecA1}, one obtains that $\phi$ is
Hadamard differentiable at
$\boldsymbol m_p\in{\cal D}^S_{0,p+1}$,
tangentially to
$\mathcal C_{0,p+1}
\subset\ell^{\infty,p+1}_0([0,1])$,
with derivative
$\phi'_{\boldsymbol m_p}$ as defined above. Further, by the Schur complement representation and the
definition of the sequential estimators, as used
in Step 4 of Section \ref{SecA1}, we have
$0\leq\widehat M_p(\lambda)
\leq\widehat M_0(\lambda)$
for all $\lambda\in[0,1]$. Thus, $\widehat{\boldsymbol M}_p
=\{(\widehat M_0(\lambda),\dots,
\widehat M_p(\lambda))^\top\}_{\lambda\in[0,1]}$
belongs to ${\cal D}^S_{0,p+1}$.
Since the limiting process
$\pmb{\mathcal I}$ in
Theorem \ref{M0, ..., Mp cvgc}
belongs to $\mathcal C_{0,p+1}$
almost surely, the functional delta method,
together with
\eqref{phi_p at g_p: Cor Mp/M0}, yields
\begin{align}\label{kue8}
\sqrt{N}
\big\{\lambda(\widehat S_p(\lambda)-S_p)
\big\}_{\lambda\in[0,1]}
\rightsquigarrow
\phi'_{\boldsymbol m_p}(\pmb{\mathcal I}).
\end{align}

Moreover, using the definition of $\pmb{\mathcal{I}}$ from Theorem \ref{M0, ..., Mp cvgc}, and equations \eqref{phi_p at g_p: Cor Mp/M0}--\eqref{eq:identity gradient fp/f0}, we obtain
\begin{align*}
    \phi'_{\boldsymbol{m}_{\!p}}(\pmb{\mathcal{I}}) &= \Big\{\frac{1}{\gamma_0} s_p\,\mathcal{M}_{p,\boldsymbol{\gamma}_p}\Sigma_p^{1/2}\,\pmb{\mathbb{B}}(\lambda)\Big\}_{\!\lambda\in[0,1]} \stackrel{d} = \big\{c \mathbb{B}(\lambda)\big\}_{\!\lambda\in[0,1]}\,,  
\end{align*}
where $\mathbb{B}$ denotes a standard Brownian motion on $[0,1]$, and 
$$
    c^2 = \frac{1}{\gamma_0^2} s_p\,\mathcal{M}_{p,\boldsymbol{\gamma}_p} \Sigma_p \,
    \mathcal{M}_{p,\boldsymbol{\gamma}_p}^\top  
    s_p^\top.
$$ 
Since $S_p = 1 - R^2_p < 1$ by assumption, and $\partial_{\gamma_0}M_p(\boldsymbol\gamma_p)\geq 1$ according to \eqref{part_der_Mp_gamma0}, we have
\[
    \mathcal{M}_{p,\boldsymbol{\gamma}_p} ^\top  s_p^\top
    \,=\,
    \big(\partial_{\gamma_0}M_p(\boldsymbol\gamma_p) - S_p, \,\partial_{\gamma_1}M_p(\boldsymbol\gamma_p), \,\dots, \,\partial_{\gamma_p}M_p(\boldsymbol\gamma_p)\big)^{\!\top} 
    \neq\, 0.
\]
As $\Sigma_p$ is positive definite and $\mathcal{M}_{p,\boldsymbol{\gamma}_p}^\top  
    s_p^\top \neq 0,$ we have $c\neq 0.$ Thus, \eqref{kue8} and an application of the continuous mapping theorem to
\begin{align*}
    \ell^{\infty}([0,1])\, \ni \boldsymbol{f} = \big\{ f(\lambda) \big\}_{\!\lambda\in[0,1]}\, \mapsto \frac{f(1)}{\int^1_0 |f(\lambda) - \lambda f(1)|\,\mathrm{d}\lambda}~,
\end{align*}
proves the weak convergence in \eqref{kue11a}. This completes the proof.\hfill\qed

\subsection{Proof of Corollary \ref{Cor:Asymptotic level test - distr. free}} This follows from Theorem \ref{piv_M_cor} and the fact that $\sqrt{N} \widehat{V}_{S_p}$ converges in distribution to an almost surely positive random variable.\hfill\qed

\subsection{Proof of Theorem \ref{thm3}}\label{SecA5}

Since $p^*$ is finite and fixed, we may assume
$N>p^*$ throughout the proof, recalling the
convention $\min\emptyset=N$ in the definition
of $\widehat p$. First, observing the decomposition 
\begin{align}\label{det83} 
    \widehat T_p(\nu) \coloneqq\frac{1- \widehat{S}_p-\nu }{\widehat{V}_{S_p}} = \frac{S_p- \widehat S_p}{\widehat{V}_{S_p}} + \frac{1- S_p- \nu }{\widehat{V}_{S_p}}
\end{align}
and the fact that the distribution of $W$ is symmetric, we obtain $\widehat T_{p^*}(\nu) \stackrel{d}{\longrightarrow}W$ if $S_{p^*} =1 -\nu,$ $\widehat T_{p^*}(\nu) \stackrel{\mathbb{P}}{\longrightarrow} \infty$ if $S_{p^*} < 1 -\nu,$ and 
$\widehat T_p(\nu) \stackrel{\mathbb{P}}{\longrightarrow} - \infty$ if $p< p^*$. Consequently, 
\[
    \mathbb{P} \big (\widehat{p} < p^*\big) = \mathbb{P}\bigg(\bigcup_{p=1}^{p^*-1} \big \{ \widehat{T}_p(\nu ) > q_\alpha(W)\big\}\bigg) \leq \sum_{p=1}^{p^*-1} \mathbb{P}\big(\widehat{T}_p (\nu ) > q_\alpha(W)\big) \longrightarrow 0,
\]
and, similarly,
\[
    \mathbb{P}\big (\widehat{p} > p^*\big) \,=\, \mathbb{P}\bigg(\bigcap_{p=1}^{p^*} \big\{\widehat{T}_p (\nu ) \leq q_\alpha(W)\big\}\bigg) \,\leq\, \mathbb{P}\big(\widehat{T}_{p^*}(\nu) \leq q_\alpha(W)\big),
\]
where the right-hand side converges to $0$ or $\alpha$ if $S_{p^*} < 1- \nu$ or $S_{p^*} = 1- \nu$, respectively. Finally, let $\alpha=\alpha_N\to0$ such that $|q_{\alpha_N}(W)|=o(\sqrt{N})$. Then, by \eqref{kue11a} and since $\widehat{V}_{S_p} = O_\mathbb{P}(N^{-1/2})$ for fixed $p$ (see Section \ref{SecA1}), we have
\[
    \widehat S_p+q_{\alpha_N}(W)\widehat V_{S_p}
    \stackrel{\Prob}{\longrightarrow} S_p\,.
\]
Thus, by the definition of $\widehat p$ and since $S_p>1-\nu$ for all $p<p^*$, we obtain $\mathbb P(\widehat p<p^*)\to0$. Similarly, if $S_{p^*}<1-\nu$, then $\mathbb P(\widehat p>p^*)\to0$. The latter conclusion also holds if $S_{p^*}=1-\nu$, since $\widehat T_{p^*}(\nu)\stackrel{d}{\longrightarrow}W$ and $q_{\alpha_N}(W)\to-\infty$. Consequently, $\mathbb P(\widehat p\neq p^*)\to0$, completing the proof. \hfill\qed

\subsection{Proof of Theorem \ref{thm4}} Under the null hypothesis, we have
$$
    \mathbb{P}_{H_0} \big (\widehat{p} > p_0 \big ) \,=\, \mathbb{P}_{p_0 \geq p^*} \Big ( \bigcap_{p=1}^{p_0} \big \{ \widehat{T}_p (\nu )\leq q_\alpha (W)  \big \} \Big ) \,\leq\, \mathbb{P}_{p_0 \geq p^*} \big( \widehat T_{p^*}(\nu)  \leq q_\alpha (W) \big ),
$$
with $\widehat T_p(\nu)$ from \eqref{det83}. By the discussion in Section \ref{SecA5}, the probability on the right-hand side converges to $\alpha$ if $S_{p^*}=1-\nu,$ and to $0$ if $S_{p^*} < 1-\nu,$ i.e.\ the test \eqref{det82} has asymptotic level $\alpha.$ Similarly, by Section \ref{SecA5}, $\widehat T_p(\nu)~{\stackrel{\mathbb{P}}{\longrightarrow}}- \infty$ for  all $p \leq p_0 < p^*.$ Consequently, under the alternative in \eqref{det92'} we obtain that 
$$
    \mathbb{P}_{H_1} (\widehat{p} > p_0) = \mathbb{P}_{p^* > p_0}\Big ( \bigcap_{p=1}^{p_0} \big \{ \widehat{T}_p (\nu )\leq q_\alpha (W)  \big \} \Big )  \longrightarrow 1 , 
$$
which proves consistency.\hfill\qed

\subsection{Proof of Theorem \ref{thm:piv_Q}}
Here, we adopt the notation from Section \ref{SecA3}.
We define the map
\begin{align*}
\phi\colon\!
\begin{cases}
\;{\cal D}^Q_{0,p+1}\rightarrow\ell^\infty([0,1]),\\
\;\boldsymbol{f}_{\!p}\mapsto
\phi(\boldsymbol{f}_{\!p})\colon\!
\begin{cases}
\;[0,1]\rightarrow\R,\\
\;\lambda\mapsto
(\phi(\boldsymbol{f}_{\!p}))(\lambda)
\coloneqq
\begin{cases}
\;\lambda f_p(\lambda)/f_{p-1}(\lambda),
& f_{p-1}(\lambda)>0,\;\lambda\in(0,1],\\
\;0, & \text{otherwise},
\end{cases}
\end{cases}
\end{cases}
\end{align*}
where
$\boldsymbol{f}_{\!p}
=\{\boldsymbol{f}_{\!p}(\lambda)\}_{\lambda\in[0,1]}
=\{(f_0(\lambda),\dots,f_p(\lambda))^\top\}_{\lambda\in[0,1]}$
and the domain is defined as
\begin{align*}
{\cal D}^Q_{0,p+1}\coloneqq
\Big\{
\boldsymbol{f}_{\!p}\in
\ell^{\infty,p+1}_0([0,1])
\,\Big|\,
0\leq f_p(\lambda)\leq f_{p-1}(\lambda)
\text{ for all }\lambda\in[0,1]
\Big\}.
\end{align*}
With $\boldsymbol{M}_{\!p}
\coloneqq(M_0,\dots,M_p)^\top\in\R^{p+1}$
and
$\boldsymbol m_p
\coloneqq\{\lambda\boldsymbol{M}_{\!p}\}_{\lambda\in[0,1]}
\in{\cal D}^Q_{0,p+1}$,
we obtain
\begin{align}\label{phi_p at g_p: Cor Mp-1/Mp}
(\phi(\boldsymbol m_p))(\lambda)
=\lambda Q_p
=\lambda\frac{M_p}{M_{p-1}},
\qquad \lambda\in[0,1],
\end{align}
where
$M_{p-1}=\det(G_{p-1})/\det(G_{p-2})>0$
by the positive definiteness of Toeplitz matrices. Further, consider the map
\begin{align*}
\widetilde{\phi}\colon
\mathcal V^Q_{p+1}\rightarrow\R,
\qquad
\boldsymbol f_{\!p}=(f_0,\dots,f_p)^\top
\mapsto\frac{f_p}{f_{p-1}},
\end{align*}
where $\mathcal V^Q_{p+1}
\coloneqq
\{
\boldsymbol f_{\!p}\in\R^{p+1}
|\,f_{p-1}>0
\}.$ The map $\widetilde{\phi}$ is continuously differentiable
on $\mathcal V^Q_{p+1},$ with gradient
\begin{align}\label{eq:identity gradient fp-1/fp}
\nabla\widetilde{\phi}
(\lambda\boldsymbol{M}_{\!p})
=\frac{1}{\lambda M_{p-1}}q_p^\top,
\qquad \lambda\in(0,1],
\end{align}
where
$q_p\coloneqq(0,\dots,0,-Q_p,1)$. Moreover, define
\begin{align*}
\phi'_{\boldsymbol m_p}\colon
\begin{cases}
\;\mathcal C_{0,p+1}
\rightarrow\ell^\infty([0,1]),\\
\;\boldsymbol h_p\mapsto
\phi'_{\boldsymbol m_p}(\boldsymbol h_p)
\colon
\begin{cases}
\;[0,1]\rightarrow\R,\\
\;\lambda\mapsto
\big(\phi'_{\boldsymbol m_p}
(\boldsymbol h_p)\big)(\lambda)
\coloneqq
\frac{1}{M_{p-1}}
q_p\boldsymbol h_p(\lambda),
\end{cases}
\end{cases}
\end{align*}
with $\mathcal C_{0,p+1}$ defined as in \eqref{setC0p+1} and equipped
with the supremum norm.

By arguments analogous to those in Sections \ref{SecA1}
and \ref{SecA3}, the map $\phi$ is Hadamard
differentiable at $\boldsymbol m_p\in{\cal D}^Q_{0,p+1}$,
tangentially to $\mathcal C_{0,p+1}
\subset\ell^{\infty,p+1}_0([0,1])$,
with derivative $\phi'_{\boldsymbol m_p}$ defined above.
Thus, the functional delta method, together with
Theorem \ref{M0, ..., Mp cvgc} and
\eqref{phi_p at g_p: Cor Mp-1/Mp}, yields
\begin{align*}
\sqrt N
\big\{\lambda(\widehat Q_p(\lambda)-Q_p)
\big\}_{\lambda\in[0,1]}
\rightsquigarrow
\phi'_{\boldsymbol m_p}(\pmb{\mathcal I})
\stackrel{d}{=}
\big\{c\mathbb B(\lambda)\big\}_{\lambda\in[0,1]},
\end{align*}
where $\mathbb B$ is a standard Brownian motion and
\[
c^2
=\frac{1}{M_{p-1}^2}
q_p\mathcal M_{p,\boldsymbol\gamma_p}
\Sigma_p
\mathcal M_{p,\boldsymbol\gamma_p}^\top
q_p^\top.
\]
To verify that $c^2>0$, note that, by \eqref{kue7},
\[
\mathcal M_{p,\boldsymbol\gamma_p}^\top q_p^\top
=
\nabla M_{p,\boldsymbol\gamma_p}
-Q_p\nabla M_{p-1,\boldsymbol\gamma_p}.
\]
The last component of this vector equals
$\partial_{\gamma_p}M_p(\boldsymbol\gamma_p)
=-2\kappa_p$ by the Schur complement
formula \eqref{Schur} and the Yule--Walker equations,
and where we assumed $\kappa_p\neq0.$ Hence, $c^2>0$ by the positive definiteness of
$\Sigma_p$ and $M_{p-1}>0$.
The claim now follows by the continuous mapping
theorem, as in Section \ref{SecA3}.\hfill\qed

\subsection{Proof of Theorem \ref{thm:PACs}}
We proceed as in Sections \ref{SecA1} and \ref{SecA3}.
Define the map
\begin{align*}
\widetilde\phi\colon
\mathbb R^{p+1}\rightarrow\mathbb R,
\qquad
\boldsymbol f_{\!p}=(f_0,\dots,f_p)^\top
\mapsto
\begin{cases}
\;\boldsymbol e_p^\top
\big((f_{|j-i|})_{i,j=0}^{p-1}\big)^{-1}
(f_1,\dots,f_p)^\top,
& \boldsymbol f_{\!p}\in\mathcal U_{p+1},\\
\;0, & \text{otherwise},
\end{cases}
\end{align*}
where $\mathcal U_{p+1}$ is defined in
Section \ref{SecA1} and $\boldsymbol e_p$
denotes the $p$th unit vector in $\mathbb R^p$.
By the Levinson--Durbin recursion,
$|\widetilde\phi(\boldsymbol f_{\!p})|\leq1$.
Moreover, $\widetilde\phi$ is continuously
differentiable on $\mathcal U_{p+1}$, with
\[
\nabla\widetilde\phi(\lambda\boldsymbol\gamma_p)
=\lambda^{-1}\boldsymbol b_p,
\qquad \lambda\in(0,1],
\]
where
$\boldsymbol b_p\coloneqq
\nabla\widetilde\phi(\boldsymbol\gamma_p)$. Next, define the map
\[
\phi\colon
\ell^{\infty,p+1}_0([0,1])
\rightarrow\ell^\infty([0,1]),
\qquad
\boldsymbol f_{\!p}\mapsto
\big\{\lambda\widetilde\phi
(\boldsymbol f_{\!p}(\lambda))
\big\}_{\lambda\in[0,1]}.
\]
With $\boldsymbol g_p$ as in
\eqref{Def hatg_p and g_p}, we have
\begin{align}\label{kue3}
(\phi(\boldsymbol g_p))(\lambda)
=\lambda\kappa_p,
\qquad \lambda\in[0,1].
\end{align}
Consider the bounded linear operator
\[
\phi'_{\boldsymbol g_p}\colon
\mathcal C_{0,p+1}\rightarrow\ell^\infty([0,1]),
\qquad
\boldsymbol h_p\mapsto
\big\{\boldsymbol b_p^\top
\boldsymbol h_p(\lambda)
\big\}_{\lambda\in[0,1]}.
\]
By arguments analogous to those in
Sections \ref{SecA1} and \ref{SecA3}, and using
$|(\phi(\boldsymbol f_{\!p}))(\lambda)|
\leq\lambda$ near zero,
$\phi$ is Hadamard differentiable at
$\boldsymbol g_p$, tangentially to
$\mathcal C_{0,p+1}
\subset\ell^{\infty,p+1}_0([0,1])$,
with derivative $\phi'_{\boldsymbol g_p}$.

Thus, the functional delta method,
together with \eqref{kue5}, \eqref{kue3},
and \eqref{eq: process cvg gammas}, yields
\begin{align*}
\sqrt N
\big\{\lambda(\widehat\kappa_p(\lambda)
-\kappa_p)\big\}_{\lambda\in[0,1]}
\rightsquigarrow
\phi'_{\boldsymbol g_p}(\pmb{\mathcal G})
\stackrel{d}{=}
\big\{c\mathbb B(\lambda)\big\}_{\lambda\in[0,1]},
\end{align*}
where $\mathbb B$ denotes a standard Brownian
motion and
$c^2=\boldsymbol b_p^\top\Sigma_p\boldsymbol b_p$.
Since the last component of
$\boldsymbol b_p$ equals
$\boldsymbol e_p^\top G_{p-1}^{-1}
\boldsymbol e_p>0$,
we have $\boldsymbol b_p\neq0$.
Hence, $c^2>0$ by the positive definiteness
of $\Sigma_p$.
The claim follows by the continuous mapping
theorem, as in Section \ref{SecA3}.
\hfill\qed

\subsection{Proof of Theorem \ref{piv_M_cor mult}}

For brevity, we only outline the main steps. Similar arguments as in Step 1 of Section
\ref{SecA1} yield
\begin{align*}
\sqrt{N}\big\{
\widehat{\boldsymbol{\nu}}_p(\lambda)
-\lambda\boldsymbol{\nu}_p
\big\}_{\lambda\in[0,1]}
\rightsquigarrow
\big\{
\boldsymbol{\Sigma}_p^{1/2}\,
\pmb{\mathbb{B}}(\lambda)
\big\}_{\lambda\in[0,1]},
\end{align*}
in $\ell^{\infty,k_{p,d}}([0,1]),$ where $\widehat{\boldsymbol{\nu}}_p(\lambda)$ is defined analogously
to $\widehat{\boldsymbol{\nu}}_p$, via the sequential
autocovariance estimates, $\boldsymbol{\Sigma}_p$ is the matrix in \eqref{det203}, and
$\pmb{\mathbb{B}}(\lambda)\coloneqq
(\mathbb{B}_1(\lambda),\dots,
\mathbb{B}_{k_{p,d}}(\lambda))^\top$
is a vector of independent standard Brownian motions on
$[0,1]$. An application of the functional delta method
then yields
\begin{align*}
    \sqrt{N}\big\{
    \lambda\big(
    \widehat{\mathscr{S}}_p(\lambda)-\mathscr{S}_p
    \big)
    \big\}_{\lambda\in[0,1]}
    \rightsquigarrow
    \big\{
    \nabla g(\boldsymbol{\nu}_p)^\top
    \boldsymbol{\Sigma}_p^{1/2}\,
    \pmb{\mathbb{B}}(\lambda)
    \big\}_{\lambda\in[0,1]}.
\end{align*}
Since $\boldsymbol{\Sigma}_p$ is non-singular and $\nabla g(\boldsymbol{\nu}_p)\neq 0$, the claim follows from the continuous mapping theorem.\hfill\qed

\end{document}